%% file: MainProfiniteOperadsPartII.tex
\documentclass[a4paper, 11pt]{article}
\usepackage{cmap}
\usepackage[english]{babel}
\usepackage[sc,osf]{mathpazo}
\usepackage[T1]{fontenc}
\usepackage[utf8]{inputenc}
\usepackage{amsfonts}       
\usepackage{amssymb}       	
\usepackage{eucal}          
\usepackage{mathrsfs}       
\usepackage{amsmath}        
\usepackage{amsthm} 
\usepackage{bbold}
\usepackage{mathtools}      
\usepackage{stmaryrd}
\usepackage{thmtools}		
\usepackage{graphicx}       
\usepackage[dvipsnames]{xcolor} 
\usepackage{accents}        
\usepackage[shortlabels]{enumitem} 
\setlist{noitemsep}
\usepackage[colorlinks=true,citecolor=PineGreen,linkcolor=BrickRed]{hyperref}  
\usepackage[capitalize]{cleveref}       
\usepackage{tikz-cd}             
\usetikzlibrary{shapes.geometric,patterns,arrows,decorations.pathreplacing}
\usepackage{caption}
\usepackage{subcaption}
\usepackage{subfiles}
\usepackage{microtype}
\usepackage{scalerel}
\usepackage{csquotes}       
\usepackage[backend=biber,style=alphabetic,sorting=nty,giveninits,urldate=long]{biblatex} 
\renewbibmacro{in:}{        
	\ifentrytype{article}{}{\printtext{\bibstring{in}\intitlepunct}}}
\AtEveryBibitem{%
	\savefield{note}{\myadd}%
	\clearfield{note}%
	\renewbibmacro{finentry}{
		\restorefield{note}{\myadd}\printfield{note}
		\finentry}%
}

\tikzset{%
	symbol/.style={%
		draw=none,
		every to/.append style={%
			edge node={node [sloped, allow upside down, auto=false]{$#1$}}}
	}
}

\tikzset{
	rot90/.style={anchor=south, rotate=90, inner sep=.5mm}
}
\tikzset{
	rot320/.style={anchor=south, rotate=320, inner sep=.5mm}
}
\tikzset{
	rot45/.style={anchor=south, rotate=45, inner sep=.5mm}
}

\addto\extrasenglish{ 
	 
	 }

\crefformat{equation}{#2(#1)#3}
\Crefformat{equation}{#2(#1)#3}

\declaretheorem[name=Theorem, numberwithin=section]{theorem}
\declaretheorem[name=Theorem]{theoremA}

\declaretheorem[name=Lemma, sibling=theorem]{lemma}
\declaretheorem[name=Proposition, sibling=theorem]{proposition}
\declaretheorem[name=Corollary, sibling=theorem]{corollary}
\declaretheorem[name=Corollary, numbered=no]{corollary*}

\declaretheorem[style=definition, name=Definition, sibling=theorem]{definition}
\declaretheorem[style=definition, name=Remark, sibling=theorem]{remark}
\declaretheorem[style=definition, name=Example, sibling=theorem]{example}

\declaretheorem[style=definition, name=Notation, sibling=theorem]{notation}

\DeclareMathOperator{\Fun}{Fun}

\DeclareMathOperator*{\colim}{colim}
\DeclareMathOperator{\id}{id}

\DeclareMathOperator{\Pro}{Pro}

\DeclareMathOperator{\Hom}{Hom}

\DeclareMathOperator{\CExp}{Cexp}
\DeclareMathOperator{\sHom}{\mathbf{hom}}
\newcommand{\sHomJoy}{\sHom_J}
\newcommand{\sHomKan}{\sHom_K}
\newcommand{\Map}{\mathrm{Map}}
\DeclareMathOperator{\Aut}{Aut}

\DeclareMathOperator{\cosk}{cosk}
\DeclareMathOperator{\sk}{sk}

\DeclareMathOperator{\SingJ}{Sing_\mathit{J}}

\newcommand{\cl}{\mathit{cl}}

\newcommand{\wh}{\widehat}
\newcommand{\wt}{\widetilde}
\newcommand{\ol}{\overline}

\newcommand{\bbN}{\mathbb N}

\newcommand{\s}{\mathbf s}
\newcommand{\dd}{\mathbf d}

\newcommand{\od}{\mathbf{od}}
\newcommand{\rod}{\mathbf{rod}}
\newcommand{\cd}{\mathbf{cd}}

\newcommand{\gd}{\mathbf{gd}}
\newcommand{\dSet}{\gd\Set}
\newcommand{\dwhSet}{\gd\wh\Set}
\newcommand{\dSpace}{\gd\Space}
\newcommand{\dProfS}{\gd\ProfS}
\newcommand{\ProfS}{\Space^{\wedge}}

\newcommand{\Cin}{\mathbf C}
\newcommand{\CinB}{\Cin_B}

\newcommand{\Din}{\mathbf D}
\newcommand{\RReedy}{\mathbf R}
\newcommand{\catE}{\mathcal E}
\newcommand{\Fin}{\mathbf{Fin}}
\newcommand{\Set}{\mathbf{Sets}}
\newcommand{\Space}{\mathbf{Spaces}}

\newcommand{\Stone}{\mathbf{Stone}}
\newcommand{\Top}{\mathbf{Top}}
\newcommand{\Op}{\mathbf{Op}}

\newcommand{\LsSet}{\mathbf{sL}}
\newcommand{\LdSet}{\mathbf{dL}}
\newcommand{\Leandspace}{\mathbf{dsL}}

\newcommand{\eqsSet}[1]{#1 \mhyphen \s\Set}
\newcommand{\eqswhSet}[1]{#1 \mhyphen \s\wh\Set}
\newcommand{\GLean}{G \mhyphen \LsSet}
\newcommand{\Omegao}{\Omega_\mathrm{o}}
\newcommand{\Omegacl}{\Omega_\mathrm{cl}}
\newcommand{\OmegaA}{\Omega_\mathrm{g}}

\newcommand{\Omegaor}{\Omega_\mathrm{or}}

\newcommand{\rodProfSf}[1]{\rod\Space^{\wedge,(#1)}}
\newcommand{\rodProfSfp}[1]{\rod\Space^{\wedge,(#1)_+}}
\newcommand{\Omegaorf}[1]{\mathbb{A}_{(#1)}}
\newcommand{\Omegaorfp}[1]{\mathbb{A}_{(#1)_+}}
\newcommand{\Mo}{M^\mathrm{o}}
\newcommand{\Lo}{L^\mathrm{o}}

\makeatletter
\newcommand{\oset}[3][0ex]{%
	\mathrel{\mathop{#3}\limits^{
			\vbox to#1{\kern-2\ex@
				\hbox{$\scriptstyle#2$}\vss}}}}
\makeatother

\newcommand{\cofarrow}{\rightarrowtail}
\newcommand{\trivcofarrow}{\oset[-.5ex]{\sim}{\rightarrowtail}}
\newcommand{\wearrow}{\oset[-.16ex]{\sim}{\longrightarrow}}

\newcommand{\fibarrow}{\twoheadrightarrow}
\newcommand{\trivfibarrow}{\oset[-.5ex]{\sim}{\twoheadrightarrow}}
\newcommand{\oface}{\oset[-0.79ex]{\circ}{\rightarrow}}
\newcommand{\iface}{\rightarrowtail}

\renewcommand{\subset}{\subseteq}

\mathchardef\mhyphen="2D

\title{Profinite completions of topological operads}
\author{Thomas Blom and Ieke Moerdijk}
\date{\today}

\begin{document}
	
	\maketitle
	\begin{abstract}
		We show that the particular profinite completion used by Boavida--Horel--Robertson \cite{Horel2017ProfiniteOperads,BoavidaHorelRobertson2019GenusZero} in their study of the Grothendieck--Teichmüller group fits in the framework of profinite completion as a left Quillen functor. More precisely, we construct a model category of profinite up-to-homotopy operads based on dendroidal objects in Quick's model category of profinite spaces and show that the construction of Boavida--Horel--Robertson extends to a left Quillen functor into this model category. We also characterize the underlying $\infty$-category of this model category and obtain a Dwyer--Kan style characterization of the weak equivalences between such profinite up-to-homotopy operads. Since this model category of profinite up-to-homotopy operads is Quillen equivalent to the one considered in \cite{BlomMoerdijk2022ProfinitePartI}, we obtain analogous results in that setting.
	\end{abstract}

\subfile{Sections/Introduction}

\subfile{Sections/Preliminaries}

\subfile{Sections/LeandSpaces}

\subfile{Sections/ProfiniteGSpaces}

\subfile{Sections/ProfinitedSpaces}

\subfile{Sections/ReedyModelStructures}

\subfile{Sections/DwyerKanEquivalences}

\subfile{Sections/UnderlyingInftyCategory}

\subfile{Sections/BoavidaHorelRobertson}

	\printbibliography
	
	\noindent{\sc Max-Planck-Institut für Mathematik, Vivatsgasse 7, 53111 Bonn, Germany.}\\
	\noindent{\emph{E-mail:} \href{mailto:blom@mpim-bonn.mpg.de}{\nolinkurl{blom@mpim-bonn.mpg.de}}}\vspace{2ex}
	
	\noindent{\sc Mathematisch Instituut, Universiteit Utrecht, Postbus 80010, 3508 TA Utrecht, The Netherlands.}\\
	\noindent{\emph{E-mail:} \href{mailto:i.moerdijk@uu.nl}{\nolinkurl{i.moerdijk@uu.nl}}}
\end{document}

%% file: Sections/Introduction.tex
\section{Introduction}

This paper is a contribution to the study of profinite completions of simplicial or topological operads. One of the difficulties in attempting to construct such a completion is that the profinite completion of spaces or simplicial sets does not commute with products. Hence a naive construction of the profinite completion by simply completing each of the spaces of operations in a given operad generally does not yield an operad. In special cases, however, the profinite completion \textit{does} preserve products up to weak homotopy equivalence. Let us call a topological operad $\mathcal{P}$ ``profinitely adapted'' if its spaces of operations have this property (see \cref{sec:Boavida-Horel-Robertson} for a precise definition). For a profinitely adapted topological operad, this naive profinite completion yields an operad up to homotopy. This is the approach taken by Horel \cite{Horel2017ProfiniteOperads} and Boavida–Horel–Robertson \cite{BoavidaHorelRobertson2019GenusZero} in their attempt to relate the profinite Grothendieck--Teichmüller group to a derived endomorphism group of an operad of braids. More explicitly, these authors formulate a universal property that a profinite completion of an (up-to-homotopy) operad should satisfy, were it to exist. They next show that for the particular operad under their consideration, the spaces of operations are profinitely adapted so the naive construction indeed yields an operad up to homotopy, which moreover satisfies their universal property. This universal property enables them to define a derived mapping space of automorphisms, which they subsequently relate to the Grothendieck–-Teichmüller group as introduced by Drinfeld \cite[\S 4]{Drinfeld1991QuasitriangularQuasiHopfAlgebras} and characterized by Fresse \cite[\S I.11]{Fresse2017HomotopyOperadsGrothendieck}, cf.\ \cite[Prop.\ 7.3]{BoavidaHorelRobertson2019GenusZero}.

The main contribution of this paper can be summarized by saying that we prove that such a profinite completion of simplicial or topological operads satisfying the universal property of Boavida, Horel and Robertson \textit{always} exists.\footnote{We show this in the context of uncoloured operads without operations of arity zero and with only the identity in arity 1, which is also the context considered in \cite{BoavidaHorelRobertson2019GenusZero}.} Moreover, this completion has very good homotopical properties; in brief, it defines a left Quillen functor from a Quillen model category of topological operads to a Quillen model category of profinite up-to-homotopy operads (``$\infty$-operads''). This shows in particular that the derived automorphism space of the operad that Boavida--Horel--Robertson define in a somewhat ad hoc way and subsequently identify with the Grothendieck--Teichmüller group fits into the general theory of derived automorphisms in Quillen model categories.

Let us now describe in more detail how we achieve these goals. Our starting point is the known fact that the category of topological operads carries a Quillen model structure, and is Quillen equivalent to the category of dendroidal sets equipped with the so-called operadic model structure, as well as to the category of dendroidal spaces equipped with the so-called complete Segal model structure \cite{CisinskiMoerdijk2013Segal}. These Quillen equivalences between topological operads, dendroidal sets and dendroidal spaces extend (in a precise sense) the well-known equivalence between the Joyal model structure on simplicial sets with quasi-categories as fibrant objects, the Bergner model structure for simplicially or topologically enriched categories \cite{Bergner2007SimplicialCats,Lurie2009HTT}, and the model category of complete Segal spaces in the sense of Rezk \cite{Rezk2001HomotopyOfHomotopy}. Complete (dendroidal) Segal spaces are a very convenient tool in many contexts, as they allow for a more direct use of the homotopy theory of spaces when working with $\infty$-categories and $\infty$-operads.

In our earlier paper “Profinite $\infty$-operads” \cite{BlomMoerdijk2022ProfinitePartI} to which this paper can be viewed as a sequel, we constructed a Quillen model category structure on the category of dendroidal profinite sets, and a Quillen adjunction between dendroidal sets and dendroidal profinite sets. This is one possible way of constructing a profinite completion functor for topological operads, now modelled by dendroidal sets.

A first step toward our goal is to construct a similar model structure on the category of dendroidal profinite spaces (rather than sets) in \cref{theorem:Main-theorem-constructing-model-structures}. For the sake of this introduction, let us denote this model category by $\Op^\wedge_\infty$. We show the following.

\begin{theoremA}\label{Theorem:A}
	The adjunction
	\[\wh{(\cdot)} : \dd\Space \rightleftarrows \Op^\wedge_\infty : U\]
	is a Quillen pair, where $\dd\Space$ is equipped with the complete Segal model structure and $U$ denotes the functor that sends a dendroidal profinite space to its underlying dendroidal space.
\end{theoremA}

This left Quillen functor is our profinite completion functor. By composing it with the Quillen equivalence between topological operads and complete dendroidal Segal spaces, we obtain a profinite completion functor for topological operads.

We then proceed to give a Dwyer–Kan style characterization of the weak equivalences in the model category $\Op^\wedge_\infty$.

\begin{theoremA}\label{Theorem:B}
	A map between fibrant objects in $\Op^{\wedge}_\infty$ is a weak equivalence precisely if it is essentially surjective and fully faithful.
\end{theoremA}

Here essential surjectivity and fully faithfulness are expressed in terms of the (profinite) spaces of colours and spaces of operations defined in \cref{sec:Dwyer-Kan}. Subsequently, we give a description of the underlying $\infty$-category of $\Op^\wedge_\infty$:

\begin{theoremA}\label{Theorem:C}
	The underlying $\infty$-category of $\Op^\wedge_\infty$ is equivalent to $\Pro_\infty(\Op_\infty^\pi)$, where $\Op_\infty^\pi$ denotes the $\infty$-category of all $\infty$-operads that are $\pi$-finite in the sense of \cref{definition:pi-finite-complete-Segal-space} and where $\Pro_\infty$ denotes the $\infty$-categorical pro-category in the sense of \cite[Def.\ 5.3.5.1]{Lurie2009HTT}.
\end{theoremA}

Given an operad whose spaces of operations are profinitely adapted, we now have two potentially distinct profinite completions: the one from \cref{Theorem:A} and the naive one, obtained by simply applying the profinite completion functor to all spaces of operations. We conclude this paper by showing that these profinite completions agree in \cref{theorem:Main-theorem-Boavida-Horel-Robertson}. We do this for a large class of topological operads, namely those that are uncoloured, have no operations in arity $0$ and have only the identity operation in arity $1$. To achieve this, we use a Postnikov-style tower of approximations to a given profinite up-to-homotopy operad. These approximations arise from a filtration in terms of the arities of the operations in the operad, inspired by the work of Göppl and Weiss \cite{Goppl2019SpectralSequenceSpaces}. We suspect that similar methods can be used in the setting of closed uncoloured operads; that is, operads with exactly one operation in arity $0$.

\paragraph{Overview of the paper}
We will start in \cref{sec:Preliminaries} by extending to dendroidal spaces some of the preliminaries about dendroidal sets from \cite{BlomMoerdijk2022ProfinitePartI}. In \cref{sec:Lean-dendroidal-spaces}, we describe the category of dendroidal profinite spaces as the pro-category of the category of so-called \emph{lean} dendroidal spaces. We characterize these by a certain finiteness property in \cref{theorem:Characterization-lean-complete-dendroidal-Segal-spaces}.

Leaving \cref{sec:Profinite-G-spaces} aside for a moment, we construct in \cref{sec:Model-structures} several model structures on categories of dendroidal profinite spaces and prove \cref{Theorem:A}. These model structures vary in two ways. First, the category of trees varies in order to also model profinite operads which are open (no nullary operations), closed (unique nullary operations) or reduced (invertible unary operations). Secondly, the weak equivalences vary: we also consider several left Bousfield localizations of the model structures we construct. We should emphasize that these Bousfield localizations are somewhat different from the usual ones in that they involve \emph{fibrantly generated} model categories rather than the more familiar cofibrantly generated ones. In particular, the localizations are not defined by inverting a small set of cofibrations. The ``complete Segal'' Bousfield localization produces the model category used in \cref{Theorem:A}, and we show that it is Quillen equivalent to the model structure for profinite $\infty$-operads considered in \cite{BlomMoerdijk2022ProfinitePartI}.

In the preceding \cref{sec:Profinite-G-spaces}, we give an alternative construction of the model structure on equivariant profinite spaces due to Quick \cite{Quick2011Continuous}. The reason for giving this construction is twofold. On the one hand, it serves as an introduction to the similar but more complicated constructions in \cref{sec:Model-structures} as just described. On the other, we obtain a new description of the fibrations of this model structure (see \cref{corollary:Fibrations-between-lean-profinite-G-spaces}). This description plays a crucial role in \cref{sec:Reedy-comparison}, where we compare the model structures constructed in \cref{sec:Model-structures} with Reedy type model structures. In this section, we also introduce a new Reedy structure on the category of trees in which all vertices have at least two incoming edges. This new Reedy structure, called the ``outer Reedy structure'', plays a central role in \cref{sec:Boavida-Horel-Robertson}.

For the complete Segal model structure on dendroidal profinite spaces, it is relatively straightforward to deduce a Dwyer--Kan style characterization of the weak equivalences (\cref{sec:Dwyer-Kan}) as well as the characterization of its underlying $\infty$-category as an $\infty$-categorical pro-category (\cref{sec:Underlying-infty}). This implies analogous results for the model structures on dendroidal profinite sets from \cite{BlomMoerdijk2022ProfinitePartI}, which are hard to prove directly in that setting.

In our final \cref{sec:Boavida-Horel-Robertson}, we address the question of how our profinite completion relates to the naive one for profinitely adapted operads considered in \cite{BoavidaHorelRobertson2019GenusZero}. We show that the two coincide for profinitely adapted operads that have no operations in arity zero and only the identity operation in arity 1. In particular, the profinite completion of \cite[\S 5]{BoavidaHorelRobertson2019GenusZero} can be extended to all such operads, whether profinitely adapted or not.

\paragraph{Acknowledgements}

Thomas Blom was supported by a grant from the Knut and Alice Wallenberg Foundation. He furthermore wishes to thank the Copenhagen Centre for Geometry and Topology (DNRF151) for its hospitality during the writing of this paper.

%% file: Sections/Preliminaries.tex
\section{Preliminaries}\label{sec:Preliminaries}

In \cite{BlomMoerdijk2022ProfinitePartI}, we were mainly concerned with categories of dendroidal sets, i.e. categories of presheaves of \emph{sets} on the category $\Omega$ of trees and its variants. In this section, we briefly recall the necessary background for presheaves of \emph{spaces} (simplicial sets) on such categories of trees.

\subsection{Categories of trees}\label{ssec:Categories-of-trees}

We refer to \cite[\S 2.1]{BlomMoerdijk2022ProfinitePartI} and the references cited there for the definition of the category $\Omega$ of trees and its full subcategories $\Omegao$ and $\Omegacl$ of open and closed trees, respectively. In this paper, we will also have use for the \emph{reduced} variant $\Omegaor$ of $\Omegao$, i.e.\ the full subcategory consisting of open trees without unary vertices. These categories of trees fit into a diagram
\begin{equation}\label{diagram:The-functors-between-categories-of-trees}
	\begin{tikzcd}
		\Omegaor \ar[r, hook, "\iota"' ,shift right]& \ar[l, shift right, "r"'] \Omegao \ar[r,hook,"o"] \ar[dr,"h"'] & \Omega \ar[d, shift right,"\cl"'] \\
		& & \Omegacl \ar[u,hook,"u"', shift right]
	\end{tikzcd}
\end{equation}
where $\iota$, $o$, and $u$ are the inclusions, $\cl$ is left adjoint to $u$, $r$ is left adjoint to $\iota$ and $h$ is the composite $\cl \circ o$. The left adjoint $r$ of $\iota$ is given by sending a tree $T$ to its ``reduction'' $r(T)$ obtained by contracting all inner edges above unary vertices. The unit $T \twoheadrightarrow r(T)$ of this adjunction is the maximal degeneracy out of $T$.

Throughout this paper, we will use the notation $\OmegaA$ to denote any of the four categories of trees defined here; the subscript $\mathrm{g}$ stands for ``generic''.

There is also a fully faithful embedding of the simplex category $\Delta$ into $\Omegao$, which sends $[n]$ to the linear tree with $n$ vertices. By postcomposing with $o \colon \Omegao \to \Omega$ or $h \colon \Omegao \to \Omegacl$, we obtain fully faithful functors $\Delta \hookrightarrow \Omega$ and $\Delta \hookrightarrow \Omegacl$. We will denote each of these three functors by $i$. Note that there is no fully faithful embedding of $\Delta$ into $\Omegaor$.

\subsection{Presheaves}\label{ssec:Presheaves-on-Omega}

A \emph{dendroidal set} is a $\Set$-valued presheaf on $\Omega$, while a \emph{dendroidal space} is a presheaf on $\Omega$ with values in simplicial sets. We will use the terminology \emph{open} dendroidal set or space, \emph{reduced open} dendroidal set or space and \emph{closed} dendroidal set or space for a presheaf on one of the subcategories mentioned above. The corresponding categories will be denoted $\dd\Set$, $\od\Set$, $\rod\Set$, $\cd\Set$, $\dd\Space$, $\od\Space$, $\rod\Space$, and $\cd\Space$. The functors between different categories of trees listed above induce functors between these categories of presheaves by restriction and Kan extension. For example, the functor $\cl \colon \Omega \to \Omegacl$ induces adjoint functors between dendroidal spaces and closed dendroidal spaces, as in
\begin{equation*}
	cl_! : \dd\Space \rightleftarrows \cd\Space : \cl^* = u_!.
\end{equation*}
The embedding $i \colon \Delta \hookrightarrow \Omega$ similarly induces an adjunction
\begin{equation}\label{equation:Embedding-sSet-into-dSet}
	i_! : \s\Set \rightleftarrows \dd\Set : i^*
\end{equation}
in which the left adjoint is fully faithful. One obtains analogous adjunctions in the case of open and closed dendroidal sets.

We will often write $\dSet$ and $\dSpace$ to denote the categories of presheaves on any of the four categories of trees defined above; i.e.\ the categories of presheaves on $\OmegaA$.

\subsection{Model structures}\label{ssec:Model-structures-preliminaries}

The category $\dd\Set$ carries the so-called \emph{operadic model structure}, whose cofibrations are normal monomorphisms and whose fibrant objects are referred to as \emph{$\infty$-operads}. This model category is Quillen equivalent to the model category of simplicial coloured operads \cite[Thm.\ 8.15]{CisinskiMoerdijk2013SimplicialOperads}. There are variations of this statement for open and closed dendroidal sets, where open dendroidal sets correspond to simplicial operads without nullary operations, while closed dendroidal sets correspond to operads with exactly one nullary operation for every colour. For details, the reader is referred to \cite{BlomMoerdijk2022ProfinitePartI} and the references cited there.

We now briefly recall several of the model structures on the category $\dd\Space$ of dendroidal spaces; for a more detailed account, the reader is referred to Chapter 12 of \cite{HeutsMoerdijk2022Trees}. The category $\dd\Space$ of dendroidal spaces carries a \emph{Reedy model structure}, based on the classical Kan--Quillen model structure on $\s\Set$ and the fact that $\Omega$ can be given the structure of a (generalized) Reedy category as in \cite[Ex.\ 10.3(3)]{HeutsMoerdijk2022Trees}. We denote this model category by $\dd\Space_R$. Its cofibrations are again the normal monomorphisms; that is, a map $X \to Y$ in $\dd\Space_R$ is a cofibration precisely if it is a monomorphism and for every tree $T \in \Omega$ and every $n \geq 0$, the group $\Aut(T)$ acts freely on the complement of the image of $X(T)_n \to Y(T)_n$. Moreover, a map $X \to Y$ in $\dd\Space_R$ is a weak equivalence if for every tree $T$, the map $X(T) \to Y(T)$ is a weak homotopy equivalence, while it is a fibration if for every tree $T$, the matching map
\begin{equation*}X(T) \to M_TX \times_{M_TY} Y(T)\end{equation*}
is a Kan fibration. For the definition of this matching map, see \cite[\S 10.2]{HeutsMoerdijk2022Trees}.\footnote{Note that \cite{HeutsMoerdijk2022Trees} uses the notation $X(\partial T)$ where we write $M_T X$.} The Reedy model structure is cofibrantly generated, with generating cofibrations given by the maps
\begin{equation}\label{equation:Generating-cofibrations-Reedy}
	\partial \Omega[T] \times \Delta[m] \cup_{\partial \Omega[T] \times \partial \Delta[m]} \Omega[T] \times \partial \Delta[m] \to \Omega[T] \times \Delta[m].
\end{equation}
Here $\Omega[T]$ denotes the dendroidal set represented by a tree $T$, viewed as a discrete dendroidal space, and $\Delta[n]$ denotes the simplicial $n$-simplex, viewed as a constant dendroidal object in $\s\Set$. In \cref{ssec:Outer-Reedy-structure}, we will study a second Reedy structure on the category of reduced open trees $\Omegaor$ of a somewhat different flavour.

The Reedy model structure on $\dd\Space$ can be further (left Bousfield) localized for the \emph{Segal} maps
\begin{equation*}
	\Omega[T] \cup_{\Omega[e]} \Omega[S] \to \Omega[S \circ_e T]
\end{equation*}
given by the grafting of a tree $T$ onto the leaf $e$ of another tree $S$. We again view the dendroidal set $\Omega[T]$ as a discrete dendroidal space. Localizing at all such Segal maps gives a model category denoted $\dd\Space_{RS}$, whose fibrant objects are called \emph{dendroidal Segal spaces}. The model category $\dd\Space_{RSC}$, whose fibrant objects are called \emph{complete dendroidal Segal spaces}, is obtained by further localizing at the map
\begin{equation*}
	\Omega[\eta] \to J.
\end{equation*}
Here $J$ denotes the nerve of the groupoid $(0 \leftrightarrow 1)$ with two objects and exactly one arrow between any ordered pair of objects, viewed as a dendroidal set via the embedding $i_! \colon \s\Set \to \dd\Set$ of \cref{equation:Embedding-sSet-into-dSet}. 

One can define the Reedy model structure and its localization at the Segal maps for the category of $\s\Set$-valued presheaves on any of the categories of trees discussed in \Cref{ssec:Categories-of-trees}. These will be denoted by $\od\Space_R$, $\od\Space_{RS}$, etc. Note that in the closed case, one should localize only for the Segal maps where $e$ is a ``very inner edge'' of $S \circ_e T$, see \cite{Moerdijk2021Closed} for details. The further localization at the map $\Omega[\eta] \to J$ can also be defined for open and closed dendroidal spaces using the embeddings $\s\Set \to \od\Set$ and $\s\Set \to \cd\Set$ of \cref{equation:Embedding-sSet-into-dSet}, defining model categories $\od\Space_{RSC}$ and $\cd\Space_{RSC}$. Note that in the closed case, one should work with the closure $\ol{\eta}$ instead of $\eta$ itself. In the reduced case there is no such embedding of $\s\Set$ into $\od\Set$, so it does not make sense to speak of the object $J$. In fact, it is explained below that completeness for reduced open dendroidal Segal spaces is in some sense automatic.

The dendroidal set $J$ is part of a cosimplicial dendroidal set $J[\bullet]$, where $J = J[1]$ and
\begin{equation*}J[n] = i_! N(0 \leftrightarrow 1 \leftrightarrow \cdots \leftrightarrow n),\end{equation*}
the nerve of the groupoid with $n+1$ objects and exactly one isomorphism between any pair of objects.
The ``singular complex'' and ``geometric realization'' with respect to $J[\bullet]$ yield a pair of adjoint functors
\begin{equation}\label{equation:Quillen-equivalence-dSpace-dSet-Sing}
	|-|_J : \dd\Space \rightleftarrows \dd\Set : \SingJ,
\end{equation}
which induces a \emph{Quillen equivalence} between $\dd\Set$ with the operadic model structure and the model category $\dd\Space_{RSC}$ (see \cite[Cor.\ 6.7 \& Prop.\ 6.11]{CisinskiMoerdijk2013Segal}). This Quillen equivalence is inverse to another one, given by simply viewing a dendroidal set as a discrete dendroidal space, respectively taking the vertices of a dendroidal space. We denote this adjunction by
\begin{equation}\label{equation:Quillen-equivalence-dSpace-dSet-0}
	d_! : \dd\Set \rightleftarrows \dd\Space : d^*.
\end{equation}
We will often omit $d_!$ from the notation and simply view $\dd\Set$ as a full subcategory of $\dd\Space$. There are similar statements for open and closed dendroidal spaces.

To obtain a variant in the reduced case, observe that $\rod\Space_{RS}$ is Quillen equivalent to the left Bousfield localization $\od\Space_{RS,\mathrm{wred}}$ of $\od\Space_{RS}$ at the map
\begin{equation*}i_!(\Delta[1]) \to i_!(\Delta[0]) = \Omega[\eta].\end{equation*}
Namely, it is shown in \cite[Thm.\ 12.62]{HeutsMoerdijk2022Trees}\footnote{See also the errata \cite{HeutsMoerdijk2022TreesErrata}.} that the adjunction
\begin{equation}\label{equation:Inclusion-reduced-Quillen-equivalence}
	\iota_! : \rod\Space_{RS} \rightleftarrows \od\Space_{RS,\mathrm{wred}} : \iota_*
\end{equation} 
induced by $\iota \colon \Omegaor \hookrightarrow \Omegao$ is a Quillen equivalence. The fibrant objects of $\od\Space_{RS,\mathrm{wred}}$ will be called \emph{weakly reduced} open dendroidal Segal spaces. Such a fibrant object is automatically complete, hence $\od\Space_{RS,\mathrm{wred}}$ is in fact a left Bousfield localization of $\od\Space_{RSC}$. Under the Quillen equivalence between $\od\Space_{RSC}$ and the model category of simplicial coloured operads, the weakly reduced open dendroidal Segal spaces correspond precisely to those simplicial coloured operads that have no nullary operations and in which all $1$-ary operations are invertible up to homotopy.

\subsection{Simplicial hom-objects}

We now discuss various simplicial hom-sets for the categories of dendroidal sets and spaces, which play an important role in the rest of this paper.

Using the tensor product of dendroidal sets \cite[Ch.\ 4]{HeutsMoerdijk2022Trees}, we define simplicial hom-sets on $\dd\Set$ by
\begin{equation*}\sHom(X,Y)_\bullet = \Hom(X \otimes i_!\Delta[\bullet], Y).\end{equation*}
Note that this is not a strict enrichment, since the tensor product of dendroidal sets is only associative up to weak equivalence. For open and closed dendroidal sets, one defines $\sHom(-,-)$ analogously.

By viewing dendroidal sets as discrete dendroidal spaces, we obtain an embedding $i_! \colon \s\Set \to \dd\Space$ and we can define a similar simplicial hom for dendroidal spaces
\begin{equation*}\sHomJoy(X,Y)_\bullet = \Hom(X \otimes i_!\Delta[\bullet], Y),\end{equation*}
where the tensor product is simply the tensor product of dendroidal sets computed degreewise. For open and closed dendroidal spaces, $\sHomJoy(-,-)$ is defined analogously.

Since $\gd\Space$ is the category of presheaves on $\OmegaA$ with values in simplicial sets, the cartesian closed structure of simplicial sets induces another simplicial hom on the category $\gd\Space$. Viewing simplicial sets as constant dendroidal spaces, this simplicial hom is given by
\begin{equation*}\sHomKan(X,Y)_\bullet = \Hom(X \times \Delta[\bullet], Y),\end{equation*}
Note that, unlike $\sHomJoy$ defined above, this is a strict simplicial enrichment of $\gd\Space$. Moreover, it is also defined on the category $\rod\Space$ of reduced open dendroidal spaces.

Both $\sHomJoy$ and $\sHomKan$ interact well with the model structures discussed above. For example, if $\OmegaA = \Omega,\Omegao,\Omegacl$, then for any normal monomorphism $X \cofarrow Y$ and any fibration $L \fibarrow K$ in $\gd\Space_{RSC}$, the pullback-power map
\begin{equation*}\sHomJoy(Y,L) \to \sHomJoy(X,L) \times_{\sHomJoy(X,K)} \sHomJoy(Y,K)\end{equation*}
is a fibration in the \emph{Joyal model structure} on $\s\Set$, which is trivial if either $X \cofarrow Y$ or $L \fibarrow K$ is a weak equivalence. The analogous statement holds for the simplicial hom-sets $\sHom(-,-)$ in the operadic model structure on $\gd\Set$.

For $\sHomKan$, the situation is even better: since the Kan--Quillen model structure on $\s\Set$ is simplicial, the Reedy model structure $\gd\Space_R$ (and hence any of its left Bousfield localizations) is simplicial. More precisely, in each of the model structures $\gd\Space_R$, $\gd\Space_{RS}$ and $\gd\Space_{RSC}$, one has that for any cofibration $X \cofarrow Y$ and any fibration $L \fibarrow K$, the map
\begin{equation*}\sHomKan(Y,L) \to \sHomKan(X,L) \times_{\sHomKan(X,K)} \sHomKan(Y,K)\end{equation*}
is a fibration in the \emph{Kan--Quillen model structure} on $\s\Set$, which is a weak homotopy equivalence if either $X \cofarrow Y$ or $L \fibarrow K$ is a weak equivalence. Note that this also holds when $\OmegaA = \Omegaor$.

Finally, for future reference, let us mention that the adjunction $d_! \dashv d^*$ gives rise to a natural isomorphism
\begin{equation}\label{equation:Simplicial-adjunction-d!-d*}
	\sHomJoy(d_!X,Y) \cong \sHom(X,d^*Y)
\end{equation}
This follows directly from the fact that $d_!$ takes a dendroidal set $X$ to the constant simplicial object with value $X$ and that the tensor product $- \otimes i_!\Delta[\bullet]$ used to define $\sHomJoy$ is the degreewise tensor product of dendroidal sets.

\subsection{Spaces of operations}\label{ssec:Spaces-of-operations}

Recall the definition of the spaces of operations of an $\infty$-operad from \cite[\S 2.1]{BlomMoerdijk2022ProfinitePartI}. Given a (closed or open) $\infty$-operad $X$, we will define its \emph{space of colours} as the maximal Kan complex $(i^*X)^{\simeq}$ contained in the quasi-category $i^*X$, where $i^*$ is as in \cref{equation:Embedding-sSet-into-dSet}.

The spaces of operations for a dendroidal Segal space are defined similarly. Let us first treat the cases of open, reduced open and general dendroidal spaces; the closed case is discussed below. In these cases, if $X$ is a dendroidal Segal space, then we call $X(\eta)$ its \emph{space of colours}. For a tuple $c_1,\ldots,c_n,d \in X(\eta)_0$ of colours, the \emph{space of operations} $X(c_1,\ldots,c_n;d)$ is defined as the pullback
\begin{equation}\label{equation:Definition-spaces-of-operations}
	\begin{tikzcd}
		X(c_1,\ldots,c_n;d) \ar[r] \ar[d] \ar[dr,very near start,phantom,"\lrcorner"] & \sHomKan(\OmegaA[C_n],X) \cong X(C_n) \ar[d,two heads] \\
		\Delta[0] \ar[r,"{(c_1,\ldots,c_n;d)}"] & \sHomKan(\partial \OmegaA[C_n],X) \cong X(\eta)^{n+1}.
	\end{tikzcd}
\end{equation}

If $X$ is a reduced open dendroidal Segal space, then we can't define the spaces of $1$-ary operations this way, since $C_1$ does not lie in $\Omegaor$. In this case, given a pair $c,d \in X(\eta)_0$, we define the space of operations $X(c;d)$ as the pullback
\begin{equation}\label{equation:Definition-spaces-of-1ary-operations-reduced}
	\begin{tikzcd}
		X(c;d) \ar[r] \ar[d] \ar[dr,very near start,phantom,"\lrcorner"] & X(\eta)^{\Delta[1]} \ar[d,two heads] \\
		\Delta[0] \ar[r,"{(c,d)}"] & X(\eta)^{2}.
	\end{tikzcd}
\end{equation}
instead.

\begin{remark}\label{remark:Spaces-of-1-ary-operations-reduced}
	To see that this definition is sensible, let $Y$ be a weakly reduced open dendroidal Segal space and suppose that $c,d \in Y(\eta)_0$ are given. Since $Y$ is local for $\Omega[C_1] \to \Omega[\eta]$, we see that $Y(\eta) \wearrow Y(C_1)$. This shows that when applying the right Quillen equivalence $\iota^* \colon \od\Space_{RS,\mathrm{wred}} \to \rod\Space_{RS}$ to $Y$, the space of operations $Y(c;d)$ defined in \cref{equation:Definition-spaces-of-operations} is weakly homotopy equivalent to the space $(\iota^*Y)(c;d)$ defined in \cref{equation:Definition-spaces-of-1ary-operations-reduced}.
\end{remark}

Recall from \cite[\S 2.3]{BlomMoerdijk2022ProfinitePartI} that for closed dendroidal sets, one has inclusions 
\begin{equation*}
	\coprod_{0\leq i \leq n} \Omega[\ol{\eta}] \hookrightarrow \partial_{\mathrm{cl}} \Omegacl[\ol{C_n}] \hookrightarrow \Omegacl[\ol{C_n}].
\end{equation*}
 In particular, for any closed dendroidal Segal space $X$, one has the following Kan fibrations
\begin{equation*}\begin{tikzcd}[column sep=scriptsize]
	X(\ol{C_n}) \cong \sHomKan(\Omegacl[\ol{C_n}], X) \ar[rr,two heads] \ar[dr,two heads,"p"'] & & \sHomKan(\partial_{\mathrm{cl}} \Omegacl[\ol{C_n}],X) \ar[dl,"q",two heads] \\
	 & X(\ol{\eta})^{n+1} &
\end{tikzcd}\end{equation*}
We call $X(\ol{\eta})$ the \emph{space of colours} of $X$. For a tuple of colours $c_1,\ldots,c_n,d \in X(\ol \eta)_0$, we define the \emph{space of operations} $X(c_1,\ldots,c_n;d)$ as the fiber of $p$ over $(c_1,\ldots,c_n,d)$ and the \emph{matching object} $X^-(c_1,\ldots,c_n;d)$ as the fiber of $q$ over $(c_1,\ldots,c_n,d)$. The induced map
\begin{equation*}X(c_1,\ldots,c_n;d) \to X^-(c_1,\ldots,c_n;d)\end{equation*}
is called the \emph{matching map}.

\begin{remark}\label{remark:Spaces-of-operations-agree}
	By \cite[Lem.\ 12.31(i)]{HeutsMoerdijk2022Trees}, the spaces of operations of a complete dendroidal Segal space $X$ are naturally equivalent to those of the dendroidal set $d^* X = X_0$, where $d^*$ is the right Quillen equivalence of \cref{equation:Quillen-equivalence-dSpace-dSet-0}. The same proof works to show that the spaces of colours of $X$ and $d^*X$ are naturally equivalent. This also applies to the open and closed case, and in the closed case it moreover follows that the matching maps of $X$ are naturally equivalent to those of $d^*X = X_0$.
\end{remark}

%% file: Sections/LeandSpaces.tex
\section{Lean dendroidal spaces}\label{sec:Lean-dendroidal-spaces}

In \cite[\S 2]{BlomMoerdijk2022ProfinitePartI} we discussed how the category of dendroidal profinite sets is itself the pro-category of a certain category of ``lean'' dendroidal sets, and we gave a homotopical characterization of these. In this section, we introduce the analogous concept of a lean dendroidal \emph{space} and derive a similar homotopical characterization.

\subsection{Skeleta and coskeleta}\label{ssec:Skeleta-and-coskeleta}
	
Recall from \cite[\S 2.2]{BlomMoerdijk2022ProfinitePartI} that a simplicial set is called \emph{lean} if it is degreewise finite and $n$-coskeletal for some $n$. Equivalently, a simplicial set $X$ is called lean if there exist a finite full subcategory $\Cin$ of $\Delta$ and a presheaf $Y \colon \Cin^{op} \to \Fin\Set$ such that $X$ is isomorphic to the right Kan extension of $Y$ along $\Cin^{op} \hookrightarrow \Delta^{op}$. The category of lean simplicial sets $\LsSet$ is a full subcategory of the category of $\Fin\Set^{\Delta^{op}}$. In particular, through the inclusion $\Fin\Set \hookrightarrow \wh\Set$ of the category of finite sets into the category of profinite sets it can be seen as a full subcategory of the category $\s\wh\Set$ of simplicial profinite sets. By \cite[Thm.\ 3.4]{BlomMoerdijk2022ProfinitePartI}, this inclusion induces an equivalence of categories
\begin{equation*}\Pro(\LsSet) \wearrow \s\wh\Set.\end{equation*}
We will generally refer to the category $\s\wh\Set$ as the category of \emph{profinite spaces}.

The definition of a lean dendroidal set is analogous: a dendroidal set $X$ is called \emph{lean} if there exist a finite full subcategory $\Cin$ of $\OmegaA$ and a presheaf $Y \colon \Cin^{op} \to \Fin\Set$ such that $X$ is isomorphic to the right Kan extension of $Y$ along the inclusion $\Cin^{op} \hookrightarrow \OmegaA^{op}$ (cf.\ \cite[Rem.\ 2.3]{BlomMoerdijk2022ProfinitePartI}). The category of lean dendroidal sets will be denoted $\LdSet$. Again, the inclusion $\LdSet \hookrightarrow \dwhSet$ induces an equivalence of categories
\begin{equation*}\Pro(\LdSet) \wearrow \dSet\end{equation*}
by \cite[Thm.\ 3.4]{BlomMoerdijk2022ProfinitePartI}. In \S 2.2 of op.\ cit., we defined the filtration 
\begin{equation*}\OmegaA^{(0)} \subset \OmegaA^{(1)} \subset \cdots \subset \OmegaA^{(n)} \subset \cdots \subset \bigcup_{n \geq 0} \OmegaA^{(n)} = \OmegaA\end{equation*}
of $\OmegaA$ by finite full subcategories, where $\OmegaA^{(n)}$ denotes the full subcategory of $\OmegaA$ spanned by those trees $T$ whose size $|T|$ is at most $n$. The \emph{size} $|T|$ of a tree $T$ is defined as the sum of the number of non-root edges and the number of vertices of $T$. By restricting and then Kan extending along the inclusion $\OmegaA^{(n)} \hookrightarrow \OmegaA$, we obtain the adjoint $n$-skeleton and $n$-coskeleton functors
\begin{equation*}\sk_n : \dSet \rightleftarrows \dSet: \cosk_n\end{equation*}
which fit into the \emph{skeletal filtration}
\begin{equation}\label{equation:skeletal-filtration}
	\sk_0(X) \hookrightarrow \sk_1(X) \hookrightarrow \cdots \hookrightarrow \sk_n(X) \hookrightarrow \cdots \hookrightarrow \colim_n \sk_nX = X
\end{equation}
and \emph{coskeleton tower}
\begin{equation}\label{equation:coskeleton-tower}
	X = \lim_n \cosk_n(X) \to \cdots \to \cosk_n(X) \to \cdots \to \cosk_1(X) \to \cosk_0(X)
\end{equation}
of a dendroidal set $X$. We warn the reader that this definition of the $n$-skeleton of a dendroidal set does not agree with the one given in \cite[\S 3.6]{HeutsMoerdijk2022Trees}.

A dendroidal set is called \emph{$n$-coskeletal} if the map $X \to \cosk_n X$ is an isomorphism. It follows immediately that $X$ is lean if and only if $X$ is degreewise finite and $n$-coskeletal for some $n$.

We will now discuss the similar property of being lean for dendroidal \emph{spaces}.

\begin{definition}\label{definition:Lean-dendroidal-space}
	A dendroidal space $X$ is called \emph{lean} if there exist a finite full subcategory $\Cin \subset \Delta \times \OmegaA$ and a presheaf $Y \colon \Cin^{op} \to \Fin\Set$ such that, as a $\Set$-valued presheaf on $\Delta \times \OmegaA$, the dendroidal space $X$ is the right Kan extension $Y$ along the inclusion $\Cin^{op} \hookrightarrow \Delta^{op} \times \OmegaA^{op}$. The full subcategory of $\dSpace$ spanned by the lean dendroidal spaces is denoted $\Leandspace$.
\end{definition}

It again follows from \cite[Thm.\ 3.4]{BlomMoerdijk2022ProfinitePartI} that the inclusion $\Leandspace \hookrightarrow \dProfS$ induces an equivalence of categories
\begin{equation*}\Pro(\Leandspace) \wearrow \dProfS.\end{equation*}

\begin{remark}\label{remark:Other-characterizations-of-lean}
	\begin{enumerate}[(a),wide,labelindent=0pt,noitemsep]
		\item\label{item1:Other-characterization-of-lean} \Cref{definition:Lean-dendroidal-space} is equivalent to saying that $X$ is degreewise finite and that there exists an $n \geq 0$ such that for every tree $T$ in $\OmegaA$, the simplicial set $X(T)$ is $n$-coskeletal and for every $p \geq 0$, the dendroidal set $X_p$ is $n$-coskeletal. Namely, in this case $X$ is the right Kan extension of its restriction to $\Delta^{op}_{\leq n} \times (\OmegaA^{(n)})^{op}$.
		\item\label{item2:Other-characterization-of-lean} In particular, if $X$ is a lean dendroidal space, then each $X(T)$ is a lean simplicial set and each $X_p$ is a lean dendroidal set. However, the converse need not be true.
	\end{enumerate}
\end{remark}

Using the Reedy structure on $\OmegaA$, one can rephrase the definition of being lean in terms of matching maps.

\begin{lemma}\label{lemma:Detecting-lean-using-matching-maps}
	Let $X$ be a dendroidal space. Then $X$ is lean if and only if 
	\begin{enumerate}[(1)]
		\item for every tree $T$, the simplicial set $X(T)$ is lean, and
		\item there exists an $n$ such that for every tree $T$ of size $|T| > n$, the matching map $X \to M_T X$ is an isomorphism.
	\end{enumerate}
\end{lemma}

\begin{proof}
	Note that for a tree $T$ of size $k+1$, the matching object $M_TX$ of a dendroidal space $X$ is defined by right Kan extending its restriction $X|_{\OmegaA^{(k)}}$ along $\OmegaA^{(k)} \hookrightarrow \OmegaA^{(k+1)}$ and subsequently evaluating at $T$. It therefore follows inductively that if $X \to M_TX$ is an isomorphism for every tree $T$ of size $|T| > n$, then $X$ is isomorphic to the right Kan extension of $X|_{\OmegaA^{(n)}}$ along $\OmegaA^{(n)} \hookrightarrow \OmegaA$. Since $\OmegaA^{(n)}$ is a finite full subcategory of $\OmegaA$, the result follows.
\end{proof}

\begin{remark}\label{remark:Detecting-lean-using-outer-matching-maps}
	The proof above uses very little about the Reedy category $\OmegaA$ itself. In fact, it only uses that the full subcategories spanned by the objects of degree at most $n$ are finite. This will be relevant in \cref{sec:Boavida-Horel-Robertson}, where we apply this lemma to a different Reedy structure on $\Omegaor$ and to some of its subcategories.
\end{remark}

We will also have use for a skeletal filtration of dendroidal spaces, combining the usual skeletal filtration for simplicial sets and the one for dendroidal sets defined above. Let us write $(\Delta \times \OmegaA)^{(n)}$ for the full subcategory of $\Delta \times \OmegaA$ spanned by those pairs $([m],T)$ where $m + |T| \leq n$. This defines a filtration of $\Delta \times \OmegaA$ by finite full subcategories. By restricting and Kan extending along the inclusion $(\Delta \times \OmegaA)^{(n)} \hookrightarrow \Delta \times \OmegaA$, we again obtain (co)skeleton functors
\begin{equation*}
	\sk_n : \dSpace \rightleftarrows \dSpace: \cosk_n
\end{equation*}
which for any dendroidal space $X$ produce a skeletal filtration and coskeleton tower as in \cref{equation:skeletal-filtration} and \cref{equation:coskeleton-tower}. Like for dendroidal sets, a dendroidal space $X$ is lean if and only if it is degreewise finite and there exists an $n$ such that $X \cong \cosk_n X$. A useful property of the $n$-skeleton is that
\begin{equation}\label{equation:Skeleton-of-representable-is-boundary}
	\sk_n(\OmegaA[T] \times \Delta[m]) \cong \partial \OmegaA[T] \times \Delta[m] \cup_{\partial \OmegaA[T] \times \partial \Delta[m]} \OmegaA[T] \times \partial \Delta[m]
\end{equation}
if $n = |T| + m -1$. In particular, the generating cofibrations \cref{equation:Generating-cofibrations-Reedy} of the Reedy model structure are precisely the maps
\begin{equation*}
	\sk_n(\OmegaA[T] \times \Delta[m]) \to \OmegaA[T] \times \Delta[m]
\end{equation*}
where $n = |T| + m - 1$.

We will also use that certain functors preserve lean objects. Recall the various functors between the categories of trees described in diagram \cref{diagram:The-functors-between-categories-of-trees} and the functors that they induce between the different categories of dendroidal spaces (cf. \cref{ssec:Presheaves-on-Omega}).

\begin{lemma}\label{lemma:Various-functors-preserving-lean-objects}
	\begin{enumerate}[(a), noitemsep]
		\item\label{item1:Functors-preserving-lean-objects} Let $X$ be a degreewise finite dendroidal space and $A$ a lean dendroidal space. Then $\sHomJoy(X,A)$ and $\sHomKan(X,A)$ are lean simplicial sets.
		\item\label{item2:Functors-preserving-lean-objects} The functor $\iota^* \colon \od\Space \to \rod\Space$ preserves lean objects.
		\item\label{item3:Functors-preserving-lean-objects} Each of the functors $o_*$, $o^*$, $o_!$, $\cl_* = u^*$, $u_*$ and $h_*$ between the categories of (closed/open) dendroidal spaces preserves lean objects.
	\end{enumerate}
\end{lemma}

\begin{proof}
	Item \ref{item1:Functors-preserving-lean-objects} follows by the same proofs as \cite[Lem.\ 2.7]{BlomMoerdijk2022ProfinitePartI} (in the case of $\sHomKan$) and \cite[Lem.\ 2.9]{BlomMoerdijk2022ProfinitePartI} (in the case of $\sHomJoy$). For \ref{item2:Functors-preserving-lean-objects}, note that $r \dashv \iota$ and hence $\iota^* \cong r_*$. The fact that $r_*$ preserves lean objects follows by the same argument as \cite[Rem.\ 2.4]{BlomMoerdijk2022ProfinitePartI}. Item \ref{item3:Functors-preserving-lean-objects} follows as in \cite[Rem.\ 2.4 \& 2.5]{BlomMoerdijk2022ProfinitePartI}.
\end{proof}

\subsection{Characterizing lean fibrant objects}

In \cite[Cor.\ 7.2.6]{BarneaHarpazHorel2017ProCategories}, a homotopical characterization of lean Kan complexes was obtained, while in \cite[\S 2.3]{BlomMoerdijk2022ProfinitePartI} a similar characterization of lean $\infty$-operads was given, where $\infty$-operads were defined as the fibrant objects in the operadic model structure on $\gd\Set$. Let us briefly recall these characterizations.

A simplicial set is called \emph{$\pi$-finite} if $\pi_0(X)$ is finite, if all homotopy groups of $X$ are finite and if there exists an $n$ such that all homotopy groups of $X$ vanish above degree $n$. In Corollary 7.2.6 of \cite{BarneaHarpazHorel2017ProCategories}, it is shown that a simplicial set $X$ is weakly equivalent to a lean Kan complex precisely if it is $\pi$-finite.

In Section 2.3 of \cite{BlomMoerdijk2022ProfinitePartI}, a similar characterization of which $\infty$-operads are weakly equivalent to lean $\infty$-operads is obtained. An (open) $\infty$-operad is called \emph{$\pi$-finite} if up to homotopy, it has finitely many colours, if all its spaces of operations are $\pi$-finite and if there exists an $n_0$ such that the spaces of operations of arity greater than $n_0$ are contractible. It is then shown in \cite[Thm.\ 2.19]{BlomMoerdijk2022ProfinitePartI} that an (open) $\infty$-operad is weakly equivalent to a lean $\infty$-operad precisely if it is $\pi$-finite.

For closed $\infty$-operads, the definition of $\pi$-finite is slightly different: instead of asking that its spaces of operations vanish above a certain arity, one asks that there exists an $n_0$ such that for all $n > n_0$, the \emph{matching maps} $X(c_1,\ldots,c_n;d) \to X^-(c_1,\ldots,c_n;d)$ are weak homotopy equivalences. Using this definition of a closed $\pi$-finite $\infty$-operad, it is shown in \cite[Thm.\ 2.33]{BlomMoerdijk2022ProfinitePartI} that a closed $\infty$-operad is weakly equivalent to a lean one precisely if it is $\pi$-finite.

Our aim is to prove an analogous result for dendroidal spaces. Recall that we defined the spaces of operations and the matching maps of a complete dendroidal Segal space in \cref{ssec:Spaces-of-operations}.

\begin{definition}\label{definition:pi-finite-complete-Segal-space}
	A fibrant object $X$ in $\dd\Space_{RSC}$, $\od\Space_{RSC}$ or $\rod\Space_{RS}$ is called \emph{$\pi$-finite} if
	\begin{enumerate}[(i),noitemsep]
		\item\label{item1:Definition-pi-finite-Segal} the set $\pi_0(X(\eta))$ is finite,
		\item\label{item2:Definition-pi-finite-Segal} for any tuple of colours $c_1,\ldots,c_n,d \in X(\eta)$, the space of operations $X(c_1,\ldots,c_n;d)$ is $\pi$-finite, and
		\item there exists an $n_0$ such that for any $n > n_0$ and any tuple of colours $c_1,\ldots,c_n,d \in X(\eta)$, the space of operations $X(c_1,\ldots,c_n;d)$ is contractible.
	\end{enumerate}
	A fibrant object $X$ in $\cd\Space_{RSC}$ is called \emph{$\pi$-finite} if it satisfies \ref{item1:Definition-pi-finite-Segal}, \ref{item2:Definition-pi-finite-Segal} and
	\begin{enumerate}
		\item[(iii')] there exists an $n_0$ such that for any $n > n_0$ and any tuple of colours $c_1,\ldots,c_n,d \in X(\eta)$, the \emph{matching map} $X(c_1,\ldots,c_n;d) \to X^-(c_1,\ldots,c_n;d)$ is a weak homotopy equivalence.
	\end{enumerate}
\end{definition}

The analogue of \cite[Thm.\ 2.19 \& 2.33]{BlomMoerdijk2022ProfinitePartI} for complete dendroidal Segal spaces will follow from the fact that the Quillen equivalences \cref{equation:Quillen-equivalence-dSpace-dSet-Sing} and \cref{equation:Quillen-equivalence-dSpace-dSet-0} interact well with the property of being lean, as expressed in the following proposition.

\begin{proposition}\label{proposition:SingJ-preserves-lean} Suppose that $\OmegaA$ is one of the categories $\Omega$, $\Omegao$ or $\Omegacl$.
	\begin{enumerate}[(a)]
		\item\label{item1:SingJ-preserves-lean} If $X$ is a lean dendroidal space, then $d^*X = X_0$ is a lean dendroidal set.
		\item\label{item2:SingJ-preserves-lean} If $A$ is a lean dendroidal set, then $\SingJ(A)$ is a lean dendroidal space.
	\end{enumerate} 
\end{proposition}

\begin{proof}
	Item \ref{item1:SingJ-preserves-lean} is obvious (cf.\ \cref{remark:Other-characterizations-of-lean}). For \ref{item2:SingJ-preserves-lean}, recall that by definition
	\begin{equation*}\SingJ(A)(T)_n = \Hom(\OmegaA[T] \otimes J[n], A).\end{equation*}
	Throughout this proof, we write $G[n]$ for the nerve of the groupoid with $n+1$ objects and exactly one isomorphism between any pair of objects; in particular, $J[n] = i_!G[n]$. Let $A$ be a lean dendroidal set and suppose that $A$ is $m$-coskeletal. Note that $\OmegaA[T] \otimes J[n] \cong \OmegaA[T] \times \mathcal{E}(G[n])$, where $\mathcal{E}$ is the functor described just above Lemma 2.8 of \cite{BlomMoerdijk2022ProfinitePartI}. By that lemma, $\mathcal{E}$ preserves degreewise finite objects, hence $\mathcal{E}(G[n])$ is a degreewise finite dendroidal set. The cartesian exponential $\CExp(\mathcal{E}(G[n]),A)$ is therefore degreewise finite and $m$-coskeletal by \cite[Lem.\ 2.7]{BlomMoerdijk2022ProfinitePartI}, and it is isomorphic to $\SingJ(A)_n$ by adjunction. By part \ref{item1:Other-characterization-of-lean} of \cref{remark:Other-characterizations-of-lean}, it remains to show that $\SingJ(A)(T)$ is an $m$-coskeletal simplicial set for any tree $T$.
	
	To see that this is the case, note that by adjunction
	\begin{equation*}
		\SingJ(A)(T)_\bullet \cong \Hom(G[\bullet],\sHom(\OmegaA[T],A)).
	\end{equation*}
	It is proved in \cite[Lem.\ 9.2]{BlomMoerdijk2020SimplicialProPublished} that for any $m$-coskeletal simplicial set $M$, the simplicial set $\Hom(G[\bullet],M)$ is again $m$-coskeletal.\footnote{Beware that in \cite[Lem.\ 9.2]{BlomMoerdijk2020SimplicialProPublished}, $J[n]$ denotes the simplicial set for which we write $G[n]$ here.} In particular, it suffices to show that $\sHom(\OmegaA[T],A)$ is $m$-coskeletal for any $m$-coskeletal dendroidal set $A$ and any tree $T$. This follows as in the proof of \cite[Lem.\ 2.9]{BlomMoerdijk2022ProfinitePartI}.
\end{proof}

\begin{theorem}\label{theorem:Characterization-lean-complete-dendroidal-Segal-spaces}
	Suppose $\OmegaA$ is one of the categories $\Omega$, $\Omegao$ or $\Omegacl$. Let $X$ be a complete dendroidal Segal space; i.e.\ a fibrant object in $\dSpace_{RSC}$. Then $X$ is weakly equivalent to a lean complete dendroidal Segal space if and only if it is $\pi$-finite.
\end{theorem}

\begin{proof}
	If $X$ is a lean complete dendroidal Segal space, then $X_0 = d^*X$ is a lean $\infty$-operad by \cref{proposition:SingJ-preserves-lean}. By \cref{remark:Spaces-of-operations-agree}, we see that $X$ is $\pi$-finite if and only if $X_0 = d^*X$ is. It follows from \cite[Thm.\ 2.19 \& 2.33]{BlomMoerdijk2022ProfinitePartI} that $X$ is $\pi$-finite. Conversely, if $X$ is $\pi$-finite, then since $X_0$ is $\pi$-finite, it follows from \cite[Thm.\ 2.19 \& 2.33]{BlomMoerdijk2022ProfinitePartI} that $X_0 \simeq A$ for some lean $\infty$-operad $A$. Then $X \simeq \SingJ(A)$, and $\SingJ(A)$ is lean by \cref{proposition:SingJ-preserves-lean}.
\end{proof}

\begin{corollary}\label{corollary:Characterizing-lean-fibrant-reduced}
	Let $X$ be a reduced open dendroidal Segal space; i.e.\ a fibrant object in $\rod\Space_{RS}$. Then $X$ is weakly equivalent to a lean reduced open dendroidal Segal space if and only if $X$ is $\pi$-finite.
\end{corollary}

\begin{proof}
	The proof that any lean reduced open dendroidal Segal space is $\pi$-finite follows as in the proof of \cite[Thm.\ 2.19]{BlomMoerdijk2022ProfinitePartI}. For the converse, let a $\pi$-finite reduced open dendroidal Segal space $X$ be given. Since $\iota^* \colon \od\Space_{RS,\mathrm{wred}} \to \rod\Space_{RS}$ is a right Quillen equivalence, there exists a weakly reduced open dendroidal Segal space $Y$ such that $\iota^*Y \simeq X$. Note that $\iota^*$ preserves the space of colours and the spaces of operations of arity greater than $1$ on the nose. Since by \cref{remark:Spaces-of-1-ary-operations-reduced} $\iota^*$ moreover preserves the spaces of $1$-ary operations up to homotopy, it follows that $Y$ must be a $\pi$-finite open complete dendroidal Segal space. By \cref{theorem:Characterization-lean-complete-dendroidal-Segal-spaces}, we see that there exists a lean weakly reduced open dendroidal Segal space $Z$ such that $Y \simeq Z$. Since $\iota^*Z$ is lean by \cref{lemma:Various-functors-preserving-lean-objects} and $X \simeq \iota^*Y \simeq \iota^*Z$, we conclude that $X$ is weakly equivalent to a lean reduced open dendroidal Segal space.
\end{proof}

%% file: Sections/ProfiniteGSpaces.tex
\section{Profinite \texorpdfstring{$G$}{G}-spaces}\label{sec:Profinite-G-spaces}
	
	In this section we give an alternative construction of a model structure originally defined by Quick \cite{Quick2011Continuous}, describing the homotopy theory of profinite $G$-spaces for a finite group $G$.\footnote{Quick actually defined such a model structure for any profinite group $G$, but in this paper we will only require the case where $G$ is finite.} The techniques used are similar to those of \cite{BlomMoerdijk2022ProfinitePartI}, so we omit most details. The main reason for presenting this alternative construction is \cref{corollary:Fibrations-between-lean-profinite-G-spaces}, which provides an explicit description of the fibrations between certain objects. This is an important ingredient in \cref{sec:Reedy-comparison}, where we compare Reedy-type model structures for dendroidal profinite spaces.
	
	Throughout this section, let $G$ be a finite group and let $\eqsSet{G}$ and $\eqswhSet{G}$ denote the categories of simplicial (profinite) sets with a \emph{right $G$-action}. We will call the objects of $\eqswhSet{G}$ \emph{profinite $G$-spaces}. It follows from \cite[\S 4]{Meyer1980Approximation} (cf.\ \cite[Thm.\ 3.3]{BlomMoerdijk2022ProfinitePartI}) that $\eqswhSet{G} \simeq \Pro(\GLean)$, where $\GLean$ denotes the category of lean simplicial sets with a right $G$-action. An object of $\GLean$ whose underlying simplicial set is a Kan complex will be called a \emph{lean $G$-Kan complex}.
	
	For $X$ and $Y$ profinite $G$-spaces, write $U_GX$ and $U_GY$ for the underlying simplicial profinite sets obtained by forgetting the $G$-action. The simplicial enrichment of $\s\wh\Set$ induces a $\eqsSet{G}$-enrichment on $\eqswhSet{G}$, given by defining $\sHom(X,Y)$ to be the simplicial set $\sHom(U_GX,U_GY)$ on which $G$ acts by conjugation. We write $\sHom_G(X,Y)$ for the $G$-fixed points of $\sHom(X,Y)$, endowing $\eqswhSet{G}$ with a simplicial enrichment.
	
	First recall that in \cite{Quick2008Profinite}, Quick defined a model structure $\s\wh\Set_Q$ on $\s\wh\Set$ describing the homotopy theory of \emph{profinite spaces}. This model structure is fibrantly generated, its cofibrations are the monomorphisms and its weak equivalences are the maps $X \to Y$ such that
		\begin{equation*}
				\sHom(Y,K) \to \sHom(X,K)
		\end{equation*}
	is a weak equivalence for any lean Kan complex $K$ (cf.\ \cite[\S 7]{BarneaHarpazHorel2017ProCategories} and \cite[Prop.\ 6.5]{BlomMoerdijk2020SimplicialProPublished}). In \cite{Quick2011Continuous}, he extended this to a model structure on $\eqswhSet{G}$ which we will denote by $\eqswhSet{G}_Q$. The cofibrations of $\eqswhSet{G}_Q$ are the \emph{normal monomorphisms}: the monomorphisms $X \to Y$ with the property that for every $n \geq 0$, the group $G$ acts freely on the complement of the image of $X_n \to Y_n$. For details on what it means for $G$ to act freely on the complement of the image of a map of profinite $G$-sets, the reader is referred to \cite[\S 3.2]{BlomMoerdijk2022ProfinitePartI}. The weak equivalences of $\eqswhSet{G}_Q$ are detected by the forgetful functor $\eqswhSet{G}_Q \to \s\wh\Set_Q$. Since any lean Kan complex can be made into a lean $G$-Kan complex and since the weak equivalences of $\eqsSet{G}$ are detected by the forgetful functor $\eqsSet{G} \to \s\Set$ to the Kan--Quillen model structure, it follows that $X \to Y$ is a weak equivalence in $\eqswhSet{G}_Q$ if and only if for every lean $G$-Kan complex $K$, the map
	\begin{equation*}\label{equation:equivariant-equivalences-detection}
		\sHom(Y,K) \to \sHom(X,K)
	\end{equation*}
	is a weak equivalence in the projective model structure on $\eqsSet{G}$. Recall that here $\sHom(-,-)$ denotes the $\eqsSet{G}$-enrichment of $\eqswhSet{G}$ defined above.

	For our alternative construction of Quick's model structure on $\eqswhSet{G}$, define the set of maps $\mathcal{Q}$ by
	\begin{equation}\label{equation:definition-Q-GSpace}
		\mathcal{Q} = \{q \colon L \to K \mid q \text{ is a trivial Kan fibration between lean $G$-Kan complexes} \}.
	\end{equation}
	Since a map in the projective model structure on $\eqsSet{G}$ is a (trivial) fibration if and only if the underlying map is a (trivial) Kan fibration, we see that this is exactly the set of trivial fibrations between lean fibrant objects in $\eqsSet{G}$. The following is now proved using the same argument as \cite[Prop.\ 4.9]{BlomMoerdijk2022ProfinitePartI}.

	\begin{proposition}\label{proposition:Normal-mono-llp-wrt-Q-Profinite-G-Spaces}
		A map in $\eqswhSet{G}$ is a normal monomorphism if and only if it has the left lifting property with respect to all maps in the set $\mathcal{Q}$.
	\end{proposition}

	\begin{proof}[Proof sketch]
		A map $A \to B$ of simplicial $G$-sets is called \emph{$n$-normal} if for each $m \leq n$, the map $A_m \to A_m$ on $m$-simplices is a monomorphism and $G$ acts freely on the complement of its image. It follows as in \cite[Lem.\ 4.8]{BlomMoerdijk2022ProfinitePartI} that up to isomorphism, every normal monomorphism in $\eqswhSet{G} \simeq \Pro(\GLean)$ admits an \emph{increasingly normal representation}: a level representation $\{f_i \colon X_i \to Y_i\}_{i \in I}$ with the property that for every $n \in \bbN$, there exists an $i \in I$ such that for any $j \leq i$, the map $f_j \colon X_j \to Y_j$ is $n$-normal. Using the well-known fact that in the projective model structure on $\eqsSet{G}$ the cofibrations are exactly the normal monomorphisms, the result now follows by the same proof as \cite[Prop.\ 4.9]{BlomMoerdijk2022ProfinitePartI}.
	\end{proof}

	Let $EG$ denote the nerve of the translation groupoid of $G$; that is, the groupoid whose objects are the elements of $G$ and with exactly one isomorphism between any two objects. Letting $G$ act by multiplication on the right, it follows that $EG$ is a degreewise finite contractible Kan complex on which $G$ acts freely. It is coskeletal since it is the nerve of a groupoid, hence it is an object of $\GLean$.
	
	\begin{lemma}
		A map $X \to Y$ in $\eqswhSet{G}$ becomes a weak equivalence in $\s\wh\Set_Q$ after forgetting the $G$-action if and only if for any lean $G$-Kan complex $K$, the induced map
		\begin{equation*}
			\sHom_G(Y \times EG, K) \to \sHom_G(X \times EG, K)
		\end{equation*}
		is a weak homotopy equivalence. 
	\end{lemma}

	\begin{proof}
		First suppose that $\sHom_G(Y \times EG, K) \to \sHom_G(X \times EG, K)$ is weak equivalence for any lean $G$-Kan complex $K$. Note that the forgetful functor $\eqsSet{G} \to \s\Set$ has a right adjoint that sends $Z$ to $\prod_{g \in G} Z$. Since this is a finite product of copies of $Z$, we see that this right adjoint takes each lean Kan complex $L$ to a lean $G$-Kan complex. In particular, we see by adjunction that for each such $L$,
		\begin{equation*}\sHom(Y \times EG, L) \to \sHom(X \times EG, L)\end{equation*}
		is a weak homotopy equivalence. We conclude that $X \times EG \to Y \times EG$, and hence $X \to Y$, is a weak equivalence underlying in $\s\wh\Set_Q$.
		
		For the converse, suppose that $X \to Y$ is a weak equivalence such that its underlying map is a weak equivalence in $\s\wh\Set_Q$ and let $K$ be a lean $G$-Kan complex. Then $\sHom(X,K)$ and $\sHom(Y,K)$ are fibrant in the projective model structure on $\eqsSet{G}$, hence
		\begin{equation*}\sHom_G(EG,\sHom(Y,K)) \to \sHom_G(EG,\sHom(Y,K))\end{equation*}
		is a weak equivalence. The result follows by adjunction.
	\end{proof}

	One can now construct a model structure on $\eqswhSet{G}$ using the strategy of \cite{BlomMoerdijk2022ProfinitePartI}. To this end, let $\mathcal{P}$ be the set
	\begin{equation*}
		\mathcal{P} = \{p \colon L \to K \mid p \text{ is a Kan fibration between lean $G$-Kan complexes} \}.
	\end{equation*}
	This set and the set $\mathcal{Q}$ defined in \cref{equation:definition-Q-GSpace} together generate Quick's model structure on $\eqswhSet{G}$.

	\begin{proposition}[\cite{Quick2011Continuous}]\label{theorem:Model-Structure-Profinite-G-Spaces-Quick}
		Let $G$ be a finite group. Then there exists a model structure on $\eqswhSet{G}$ in which a map is a cofibration if and only if it is a normal monomorphism and a map is a weak equivalence if and only if its underlying map in $\s\wh\Set_Q$ is a weak equivalence.
	\end{proposition}

	\begin{proof}[Proof sketch]
		The proof of \cite[Thm.\ 5.12]{BlomMoerdijk2022ProfinitePartI} also applies here, taking $\mathcal{P}$ and $\mathcal{Q}$ the sets defined above and using $EG$ instead of the object denoted by $E$ in loc.\ cit.
	\end{proof}
	
	The following technical result is now a direct consequence of our construction.
	
	\begin{corollary}\label{corollary:Fibrations-between-lean-profinite-G-spaces}
		Let $G$ be a finite group. If $f \colon X \to Y$ is a map between lean $G$-Kan complexes, then $f$ is a (trivial) fibration in Quick's model structure on $\eqswhSet{G}$ if and only if its underlying map in $\s\Set$ is a (trivial) Kan fibration.
	\end{corollary}

	\begin{proof}
		Let $f$ be a map between lean $G$-Kan complexes. Then $f$ is a Kan fibration if and only if $f$ lies in $\mathcal{P}$, while it is a trivial Kan fibration if and only if it lies in $\mathcal{Q}$. Since by the proof of \cref{theorem:Model-Structure-Profinite-G-Spaces-Quick}, the sets $\mathcal{P}$ and $\mathcal{Q}$ serve as generating (trivial) fibrations of Quick's model structure on $\eqswhSet{G}$, the result follows.
	\end{proof}

	The following lemma, which follows by the same proof as \cite[Prop.\ 5.18]{BlomMoerdijk2022ProfinitePartI}, will be used in \cref{sec:Boavida-Horel-Robertson}.
	
	\begin{lemma}\label{lemma:Profinite-G-Space-Filtered-colimit-lean-G-Kan}
		Let $G$ be a finite group and let $X$ be fibrant in Quick's model structure on $\eqswhSet{G}$. Then $X$ is a cofiltered limit of lean $G$-Kan complexes. \qed
	\end{lemma}

%% file: Sections/ProfinitedSpaces.tex
\section{Model structures for profinite dendroidal spaces}\label{sec:Model-structures}

In this section, we will construct various model structures on the category $\dd\ProfS$ of profinite dendroidal spaces and its variants, and prove some of their basic properties. Since the construction of these model structures is virtually identical for each of the indexing categories $\Omega$, $\Omegao$, $\Omegacl$ and $\Omegaor$, we again use $\OmegaA$ to denote any of these categories and write $\dSpace$ and $\dProfS$ for the category of $\s\Set$-valued and $\s\wh\Set$-valued presheaves on $\OmegaA$, respectively. Many of the proofs given in \cite{BlomMoerdijk2022ProfinitePartI} carry over to the current setting without much change, so throughout this section the reader is often referred to op.\ cit.\ for details. The arguments are similar to those of \cref{sec:Profinite-G-spaces}, though slightly more involved since there is no analogue of the object $EG$; that is, there is no \emph{lean} object in $\dSpace$ that is normal, contractible and fibrant. However, a degreewise finite object having these properties still exists, and suffices for our purposes (cf.\ \cref{lemma:Degreewise-finite-normalization-point}).
\subsection{Normal monomorphisms}

The cofibrations in each of the model structures on $\dProfS$ will be the \emph{normal monomorphisms}.

\begin{definition}\label{definition:Normal-monomorphism-dSpaces}
	\begin{enumerate}[(a)]
		\item A morphism $X \to Y$ in $\dSpace$ or $\dProfS$ is called a \emph{normal monomorphism} if for each tree $T$, the map $X(T) \to Y(T)$ is a normal monomorphism in $\eqsSet{\Aut(T)}$ or $\eqswhSet{\Aut(T)}$.
		\item A (profinite) dendroidal space $X$ is called \emph{normal} if $\varnothing \to X$ is a normal monomorphism.
	\end{enumerate}
\end{definition}

The proof of \cite[Lem.\ 12.1]{HeutsMoerdijk2022Trees} goes through for each variation of $\OmegaA$ considered here, showing that the cofibrations in the Reedy model structure on $\dSpace$ are precisely the normal monomorphisms. In particular, for each model structure on $\dSpace$ considered in \cref{ssec:Model-structures-preliminaries}, the cofibrations agree with the normal monomorphisms.

We can now mimic \cite[\S 4]{BlomMoerdijk2022ProfinitePartI} to obtain alternative characterizations of the normal monomorphisms in $\dProfS$. Recall the definition of $\sk_n$ from \cref{ssec:Skeleta-and-coskeleta}.

\begin{definition}
	\begin{enumerate}[(a)]
		\item A morphism $f \colon X \to Y$ of dendroidal spaces is called \emph{$n$-normal} if $\sk_n X \to \sk_n Y$ is a normal monomorphism.
		\item Let $I$ be a codirected poset and let $\{f_i \colon X_i \to Y_i\}_{i \in I}$ be a natural map between diagrams of lean dendroidal spaces indexed by $I$. Then $\{f_i \colon X_i \to Y_i\}_{i \in I}$ is \emph{increasingly normal} if for any $n \in \bbN$, there exists an $i \in I$ such that for any $j \leq i$, the map $f_j \colon X_j \to Y_j$ is $n$-normal. A map in $\Pro(\Leandspace)$ is said to \emph{admit an increasingly normal representation} if it can be represented by a natural transformation of this kind.
	\end{enumerate}
\end{definition}

\begin{remark}
	By the same argument as \cite[Lem.\ 4.6]{BlomMoerdijk2022ProfinitePartI}, it follows that $X \to Y$ is $n$-normal if and only if for each tree $T$ and each natural number $m$ such that $|T| + m \leq n$, the map $X(T)_m \to Y(T)_m$ is a monomorphism of sets and $\Aut(T)$ acts freely on the complement of its image.
\end{remark}

The following is proved exactly as in \cite[Lem.\ 4.8]{BlomMoerdijk2022ProfinitePartI}.

\begin{proposition}
	Let $f \colon X \to Y$ be a map of profinite dendroidal spaces. Then $f$ is a normal monomorphism if and only if, up to isomorphism, it admits an increasingly normal representation. \qed
\end{proposition}

Now view $\dSpace$ as endowed with any of the model structures of \cref{ssec:Model-structures-preliminaries}, and define the set of maps
\begin{equation}\label{equation:definition-Q-dSpace}
	\mathcal{Q} = \{q \colon B \to A \mid q \text{ is a trivial fibration between lean fibrant objects in } \dSpace \}.
\end{equation}

Just like in \cref{proposition:Normal-mono-llp-wrt-Q-Profinite-G-Spaces}, the following characterization of the normal monomorphisms in $\dProfS$ follows by the same proof as \cite[Lem.\ 4.9]{BlomMoerdijk2022ProfinitePartI}.

\begin{proposition}
	A morphism $X \to Y$ of profinite dendroidal spaces is a normal monomorphism if and only if it has the left lifting property with respect to all maps in $\mathcal{Q}$. \qed
\end{proposition}

\subsection{Weak equivalences}

Throughout this section, let $\dSpace$ be endowed with any of the model structures of \cref{ssec:Model-structures-preliminaries}. Unlike in the case of $\eqsSet{G}$, there need not be a lean normal object $E$ in $\dSpace$ for which the map to the terminal object is a trivial fibration (cf.\ \cite[Ex.\ 4.3]{BlomMoerdijk2022ProfinitePartI}). However, one can construct a degreewise finite one, where degreewise finite means that $E(T)_n$ is finite for every tree $T$ and every $n \geq 0$.

\begin{lemma}\label{lemma:Degreewise-finite-normalization-point}
	There exists a degreewise finite normal dendroidal space $E$ such that the map $E \to *$ to the terminal object is a trivial fibration in $\dSpace$.
\end{lemma}

\begin{proof}
	Recall that the generating cofibrations of $\dd\Space$ are the maps of the form
	\begin{equation}\label{equation:Generating-cofibrations-Reedy-2}
		\partial \OmegaA[T] \times \Delta[m] \cup_{\partial \OmegaA[T] \times \partial \Delta[m]} \OmegaA[T] \times \partial \Delta[m] \to \OmegaA[T] \times \Delta[m].
	\end{equation}
	As in \cite[Lem.\ 5.1]{BlomMoerdijk2022ProfinitePartI}, the result follows by modifying Quillen's small object argument where one only attaches maps of the kind \cref{equation:Generating-cofibrations-Reedy-2} with $|T| + m = n$ at the $n$-th stage. Note that this uses the fact that
	\begin{equation*}\sk_n (\OmegaA[T] \times \Delta[m]) \cong \partial \OmegaA[T] \times \Delta[m] \cup_{\partial \OmegaA[T] \times \partial \Delta[m]} \OmegaA[T] \times \partial \Delta[m] \end{equation*}
	if $n = |T| + m -1$ (see \cref{equation:Skeleton-of-representable-is-boundary}).
\end{proof}

For $X = \{X_i\}$ and $Y = \{Y_j\}$ in $\Pro(\Leandspace) \simeq \dSpace$, we define simplicial hom-sets by
\begin{equation*}
	\sHomJoy(X,Y) = \lim_j \colim_i \sHomJoy(X_i,Y_j) \quad\text{and}\quad \sHomKan(X,Y) = \lim_j \colim_i \sHomKan(X_i,Y_j).
\end{equation*}
Using \cref{lemma:Various-functors-preserving-lean-objects}, one can prove by a similar argument as \cite[Lem.\ 5.4]{BlomMoerdijk2022ProfinitePartI} that for any normal $X$ in $\dProfS$ and any lean fibrant object $A$ in $\dSpace$, the simplicial set $\sHomKan(X,A)$ is a Kan complex. When using the (complete) Segal model structure, the same argument moreover shows that $\sHomJoy(X,A)$ is an $\infty$-category.

In the following, fix a degreewise finite normalization $E \trivfibarrow *$ in the sense of \cref{lemma:Degreewise-finite-normalization-point}. Since $E$ is degreewise finite, we can view it as an object of $\dProfS$ through the inclusion $\Fin\Set \hookrightarrow \wh\Set$.

\begin{definition}\label{definition:profinite-weak-equivalences}
	A map $X \to Y$ in $\dProfS$ is called a \emph{weak equivalence} if for any lean fibrant object $A$ in $\dSpace$, the induced map
	\begin{equation*}
		\sHomKan(Y \times E, A) \to \sHomKan(X \times E, A)
	\end{equation*}
	is a weak equivalence of Kan complexes.
\end{definition}

It turns out that this definition does not depend on the choice of $E$ (cf.\ \cref{remark:WEs-independent-of-normalization}). However, until then we fix a choice for the degreewise finite normalization $E$.

\subsection{Construction of the model structures}

Let $\dSpace$ be equipped with one of the model structures of \cref{ssec:Model-structures-preliminaries}. In analogy with the definition of $\mathcal{Q}$ in \Cref{equation:definition-Q-dSpace} above, we define the following set of maps
\begin{equation}\label{equation:definition-P-dSpace}
	\mathcal{P} = \{p \colon B \to A \mid p \text{ is a fibration between lean fibrant objects in } \dSpace \}.
\end{equation}
These sets together generate the following model structure on $\dProfS$, which can be seen as a profinite analogue of the one on $\dSpace$. The proof of \cite[Thm.\ 5.12]{BlomMoerdijk2022ProfinitePartI} goes through without a change in this setting.

\begin{theorem}\label{theorem:Main-theorem-constructing-model-structures}
	There exists a fibrantly generated model structure on $\dProfS$ in which the cofibrations are the normal monomorphisms and the weak equivalences are as defined in \cref{definition:profinite-weak-equivalences}. The sets $\mathcal{P}$ and $\mathcal{Q}$ serve as generating fibrations and generating (trivial) fibrations, respectively, and this model structure is left proper. \qed
\end{theorem}

\begin{notation}\label{notation:Subscripts-dProfS}
	Recall from \cref{ssec:Model-structures-preliminaries} that we use subscripts to denote which model structure on the category of dendroidal spaces we are considering. The model structures obtained from \Cref{theorem:Main-theorem-constructing-model-structures} will be denoted using the same subscripts. For example, $\od\ProfS_{RS}$ denotes the model structure on the category of open profinite dendroidal spaces obtained by applying \Cref{theorem:Main-theorem-constructing-model-structures} to $\od\Space_{RS}$.
\end{notation}

\begin{remark}\label{remark:WEs-independent-of-normalization}
	Note that the sets of generating (trivial) fibrations $\mathcal{P}$ and $\mathcal{Q}$ are independent of the choice of degreewise finite normalization $E \trivfibarrow *$. Since the fibrations and trivial fibrations completely determine a model structure, it follows that the weak equivalences of \cref{definition:profinite-weak-equivalences} must be independent of the choice of $E$.
\end{remark}

\subsection{Some basic properties}\label{ssec:Basic-properties}

We will now discuss some properties of the model structures constructed above. Note that the first few are analogues of results from \cite[\S 5]{BlomMoerdijk2022ProfinitePartI} for presheaves with values in profinite sets.

In the propositions below, $\dProfS$ denotes any of the model structures obtained by applying \cref{theorem:Main-theorem-constructing-model-structures} to any of the model structures on $\dSpace$ mentioned in \cref{ssec:Model-structures-preliminaries}. The following propositions follow by the same proofs as Propositions 5.9, 5.14 and 5.18 of \cite{BlomMoerdijk2022ProfinitePartI}.

\begin{proposition}\label{proposition:Weak-equivalences-stable-under-cofiltered-limits}
	Weak equivalences in $\dProfS$ are stable under cofiltered limits. \qed
\end{proposition}

\begin{proposition}\label{proposition:Simplicial-model-structures}
	The simplicial enrichment $\sHomKan$ makes $\dProfS$ into a simplicial model category with respect to the Kan--Quillen model structure on $\s\Set$. \qed
\end{proposition}

\begin{proposition}\label{proposition:Fibrant-objects-are-pro-fibrant}
	Let $X$ be a fibrant object in $\dProfS$. Then $X$ is a cofiltered limit of lean objects that are fibrant in $\dSpace$. \qed
\end{proposition}

We now give an alternative characterization of the weak equivalences between cofibrant objects, independent of $E$.

\begin{proposition}
	Let $f \colon X \to Y$ be a map between cofibrant objects in $\dProfS$. Then $f$ is a weak equivalence if and only if for every lean fibrant object $A$ in $\dSpace$, the map
	\begin{equation*}f^* \colon \sHomKan(Y, A) \to \sHomKan(X, A)\end{equation*}
	is a weak equivalence of Kan complexes.
\end{proposition}

\begin{proof}
	Since $X \times E \to X$ and $Y \times E \to Y$ are weak equivalences between cofibrant objects, it follows that $\sHomKan(Y \times E, A) \to \sHom(X \times E, A)$ is a weak equivalence if and only if $\sHomKan(Y, A) \to \sHomKan(X,A)$ is.
\end{proof}

The following is shown by the same argument and Corollary 5.17 of \cite{BlomMoerdijk2022ProfinitePartI}.

\begin{proposition}\label{proposition:Profinite-completion-left-Quillen}
	The adjunction
	\begin{equation*} \wh{(\cdot)} : \dSpace \rightleftarrows \dProfS : U \end{equation*}
	given by the profinite completion functor and its right adjoint is a Quillen pair. \qed
\end{proposition}

\begin{definition}\label{definition:Profinite-completion-operad}
	Let $X$ be an $\infty$-operad, modelled as an (open/closed) complete dendroidal Segal space. Then its \emph{profinite completion} is defined as a fibrant replacement of $\wh X$ in $\dProfS_{RSC}$. If $X$ is modelled by a reduced open dendroidal Segal space, then its \emph{profinite completion} is defined as a fibrant replacement of $\wh X$ in $\rod\ProfS_{RS}$.
\end{definition}

To make the previous definition homotopy-invariant, one might wish to replace $X$ with a normalization before applying $\wh{(\cdot)}$. The following lemma ensures that this is not necessary, since the functor $\wh{(\cdot)}$ is already homotopy-invariant.

\begin{lemma}\label{lemma:Profinite-completion-preserves-all-equivalences}
	The profinite completion functor $\wh{(\cdot)}$ preserves all weak equivalences. In particular, if $X \to Y$ is a weak equivalence in $\dSpace$ between degreewise finite objects, then it is also a weak equivalence in $\dProfS$.
\end{lemma}

\begin{proof}
	To see that $\wh{(\cdot)}$ preserves all weak equivalences, it suffices to show that it takes $X \times E \trivfibarrow X$ to a weak equivalence. Since $E$ is degreewise finite, it follows that $\wh{X \times E} \cong \wh{X} \times E$, hence this is indeed the case. The second statement now follows since
	\begin{equation*}\begin{tikzcd}
		\Fin\Set^{\Delta^{op} \times \OmegaA^{op}} \ar[dr,hook] \ar[d,hook] & \\
		\Set^{\Delta^{op} \times \OmegaA^{op}} \ar[r,"\wh{(\cdot)}"'] & \wh\Set^{\Delta^{op} \times \OmegaA^{op}}
	\end{tikzcd}\end{equation*}
	commutes.
\end{proof}

We will now restrict our attention to the model structures for complete profinite dendroidal Segal spaces. Let us first note that $\sHomJoy(-,-)$ also behaves well with respect to these model structures.

\begin{proposition}\label{proposition:Joyal-hom-detects-weak-equivalences}
	Let $\OmegaA$ be one of the categories $\Omega$, $\Omegao$ or $\Omegacl$, and let $E \trivfibarrow *$ be any degreewise finite normalization in the sense of \cref{lemma:Degreewise-finite-normalization-point}. Then a map $f \colon X \to Y$ in $\dProfS_{RSC}$ is a weak equivalence if and only if for every lean fibrant object $A$ in $\dSpace_{RSC}$, the map
	\begin{equation*}\sHomJoy(Y \times E, A) \to \sHomJoy(X \times E, A)\end{equation*}
	is a weak equivalence of $\infty$-categories. Moreover, for any cofibration $Z \cofarrow W$ and any fibration $L \fibarrow K$ in $\dSpace_{RSC}$, the map
	\begin{equation*}\sHomJoy(Y,L) \to \sHomJoy(X,L) \times_{\sHomJoy(X,K)} \sHomJoy(Y,K)\end{equation*}
	is a fibration in Joyal's model structure $\s\Set_J$, which is trivial if either $Z \cofarrow W$ or $L \fibarrow K$ is.
\end{proposition}

\begin{proof}
	The proof of \cref{theorem:Main-theorem-constructing-model-structures} to construct the model structure $\dProfS_{RSC}$ goes through without a change if one defines the weak equivalences to be those maps $f \colon X \to Y$ for which
	\begin{equation}\label{equation:profinite-weak-equivalence-Joyal-version}
		\sHomJoy(Y \times E, A) \to \sHomJoy(X \times E, A)
	\end{equation}
	is a weak equivalence of $\infty$-categories for every lean fibrant $A$ in $\dSpace_{RSC}$, using $\sHomJoy$ in place of $\sHomKan$ everywhere in the proof. The resulting model structure has the same generating (trivial) fibrations as $\dProfS_{RSC}$, hence it must be the same model structure. In particular, a map is a weak equivalence in $\dProfS_{RSC}$ if and only if \cref{equation:profinite-weak-equivalence-Joyal-version} is a weak equivalence of $\infty$-categories for every lean fibrant $A$. The second statement now follows by the same proof as \cite[Prop.\ 5.14]{BlomMoerdijk2022ProfinitePartI}.
\end{proof}

The model structures constructed here are related through several Quillen pairs. For example, recall the functors $o$, $u$, $h$ and $\iota$ from diagram \cref{diagram:The-functors-between-categories-of-trees}. Using \cref{lemma:Various-functors-preserving-lean-objects} one obtains right Quillen functors
\begin{align*}
	&o^* \colon \dd\ProfS_{RSC} \to \od\ProfS_{RSC}, \\ 
	&u^* \colon \dd\ProfS_{RSC} \to \cd\ProfS_{RSC} \quad \text{and} \\
	&h_* \colon \od\ProfS_{RSC} \to \cd\ProfS_{RSC}
\end{align*}
by the same arguments as \cite[Prop.\ 5.21]{BlomMoerdijk2022ProfinitePartI}. Since the right Quillen functor
\begin{equation*}\iota^* \colon \od\Space_{RS,\mathrm{wred}} \to \rod\Space_{RS}\end{equation*}
preserves lean objects, it induces a similar Quillen pair. We will show that the latter is in fact a Quillen equivalence, allowing us to identify reduced open profinite dendroidal Segal spaces with a left Bousfield localization of the model category of complete open profinite dendroidal Segal spaces.

\begin{theorem}\label{theorem:Reduced-profinite-equivalent-to-weakly-reduced-profinite}
	The inclusion $\iota \colon \Omegaor \hookrightarrow \Omegao$ induces a Quillen equivalence
	\begin{equation*}\iota_! : \rod\ProfS_{RS} \rightleftarrows \od\ProfS_{RS,\mathrm{wred}} : \iota^*.\end{equation*}
\end{theorem}

\begin{proof}
	It suffices to show that the derived counit is an equivalence and that $\iota_!$ detects weak equivalences between cofibrant objects. Note that since $\iota^*$ preserves lean objects, $\iota_!$ must preserve cofiltered limits.
	
	To see that the derived counit is an equivalence, let $X$ be fibrant in $\od\ProfS_{RS,\mathrm{wred}}$. Since $\iota^*$ preserves normal monomorphisms, the derived counit is weakly equivalent to the map $\iota_! \iota^*(X \times E) \to X \times E$.
	By \cref{proposition:Fibrant-objects-are-pro-fibrant}, we may write $X = \lim_i X_i$ as an inverse limit of lean objects that are fibrant in $\rod\ProfS_{RS,\mathrm{wred}}$.
	Observe that for any lean fibrant $A$ in $\od\Space_{RS,\mathrm{wred}}$, the map
	\begin{equation*}\sHomKan(X_i \times E,A) \to \sHomKan(\iota_!\iota^*(X_i \times E), A) \cong \sHomKan(\iota^*(X_i \times E), \iota^* A)\end{equation*}
	is a weak equivalence because of the Quillen equivalence \cref{equation:Inclusion-reduced-Quillen-equivalence}, hence $\iota_! \iota^* (X_i \times E) \to X_i \times E$ is a weak equivalence. By \cref{proposition:Weak-equivalences-stable-under-cofiltered-limits}, we conclude that the derived counit
	\begin{equation*}\iota_! \iota^*( X \times E) \cong \lim_i \iota_! \iota^*(X_i \times E) \to \lim_i X_i \times E \cong X \times E\end{equation*}
	is a weak equivalence.
	
	Finally, to show that $\iota_!$ detects weak equivalences between cofibrant objects, suppose that $X \to Y$ is a map between cofibrant objects in $\rod\ProfS_{RS}$ such that $\iota_!X \to \iota_!Y$ is a weak equivalence in $\od\ProfS_{RS,\mathrm{wred}}$. 
	Let $A$ be lean and fibrant in $\rod\Space_{RS}$.
	It follows as in the proof of \cref{corollary:Characterizing-lean-fibrant-reduced} that there exists a lean fibrant object $B$ in $\od\Space_{RS,\mathrm{wred}}$ such that $A \simeq \iota^*B$ in $\rod\Space_{RS}$.
	Then $A \simeq \iota^*B$ in $\rod\ProfS_{RS}$ by \cref{lemma:Profinite-completion-preserves-all-equivalences}, hence $\sHomKan(Y,A) \to \sHomKan(X,A)$ is a weak equivalence if and only if
	\begin{equation*}\sHomKan(Y,\iota^*B) \cong \sHomKan(\iota_!Y,B) \to \sHomKan(\iota_!X,B) \cong \sHomKan(X,\iota^*B)\end{equation*}
	is a weak equivalence.
	The latter holds by the assumption, so we conclude that $X \to Y$ is a weak equivalence.
\end{proof}

We will now show the equivalence between $\dProfS_{RSC}$ and the model structure for profinite $\infty$-operads on $\gd\wh\Set$ (cf.\ \cite[Thm.\ 5.12]{BlomMoerdijk2022ProfinitePartI}).
This turns out to be a consequence of \cref{proposition:SingJ-preserves-lean} and the way in which these model structures are constructed.
Recall that $\LdSet$ and $\Leandspace$ denote the categories of lean dendroidal sets and lean dendroidal spaces, respectively.
By \cref{proposition:SingJ-preserves-lean}, we have functors
\begin{equation*} \SingJ \colon \LdSet \to \Leandspace \quad \text{and} \quad d^* \colon \Leandspace \to \LdSet. \end{equation*}
Since these preserve finite limits, they induce adjunctions
\begin{equation*} |-|_J : \dd\ProfS \rightleftarrows \dd\wh\Set : \SingJ \quad \text{and} \quad d_! : \dd\wh\Set \rightleftarrows \dd\ProfS : d^*. \end{equation*}
One obtains analogous adjunctions in the open and closed cases.

\begin{theorem}\label{theorem:Profinite-Segal-spaces-vs-profinite-operads}
	Let $\OmegaA$ be one of the categories $\Omega$, $\Omegao$ or $\Omegacl$. Then both adjunctions
	\begin{equation*}|-|_J : \dProfS_{RSC} \rightleftarrows \gd\wh\Set : \SingJ \quad \text{and} \quad d_! : \gd\wh\Set \rightleftarrows \dProfS_{RSC} : d^*\end{equation*}
	are Quillen equivalences, where $\gd\wh\Set$ is endowed with the operadic model structure of \cite[Thm.\ 5.12]{BlomMoerdijk2022ProfinitePartI}.
\end{theorem}

\begin{proof}
	Since $\SingJ \colon \gd\Set \to \gd\Space_{RSC}$ and $d^* \colon \gd\Space_{RSC} \to \gd\Set$ are right Quillen, \cref{proposition:SingJ-preserves-lean} implies that they must preserve the generating (trivial) fibrations of $\gd\wh\Set$ and $\dProfS_{RSC}$.
	In particular, the adjunctions in the statement of this theorem are Quillen pairs.
	Since $d^* \circ \SingJ$ is naturally isomorphic to the identity, the theorem follows if we can show that $d_! \dashv d^*$ is a Quillen equivalence.
	
	The natural isomorphism $\sHomJoy(d_!X,Y) \cong \sHom(X,d^*Y)$ from \cref{equation:Simplicial-adjunction-d!-d*} induces a similar natural isomorphism in the profinite case, where $\sHom$ denotes the simplicial hom of the category of dendroidal (profinite) sets. The proof that 
	\begin{equation*}d_! : \gd\wh\Set \rightleftarrows \dProfS_{RSC} : d^*\end{equation*}
	is a Quillen equivalence now proceeds the same as the proof of \cref{theorem:Reduced-profinite-equivalent-to-weakly-reduced-profinite}: Replacing $\sHomKan$ by $\sHomJoy$ and $\sHom$ as necessary and noting that $d^*$ preserves normal objects, it follows that the derived counit is an equivalence. In the proof that $d_!$ detects weak equivalences between cofibrant objects, to see that any lean fibrant object $A$ in $\gd\Set$ is equivalent to an object of the form $d^*B$ with $B$ lean and fibrant in $\dSpace$, simply note that $A \cong d^*\SingJ(A)$.
\end{proof}

\begin{remark}\label{remark:Profinite-completions-agree}
	The functor $d^* \colon \dProfS \to \dwhSet$ is given by sending a simplicial object $\Delta^{op} \to \dwhSet$ to its zeroth object, while $d_! \colon \dwhSet \to \dProfS$ sends a dendroidal profinite set $X$ to the constant dendroidal profinite space with value $X$. In particular, both $d^*$ and $d_!$ commute with the profinite completion functors. This shows that the profinite completions $\gd\Space_{RSC} \to \dProfS_{RSC}$ and $\dd\Set \to \dwhSet$ agree under the Quillen equivalence of \cref{theorem:Profinite-Segal-spaces-vs-profinite-operads}.
\end{remark}

%% file: Sections/ReedyModelStructures.tex
\section{Profinite Reedy model structures}\label{sec:Reedy-comparison}

The category $\Omega$ of trees and its variants are well-known examples of (generalized) Reedy categories. In particular, presheaves on these categories with values in a fibrantly or cofibrantly generated model category admit a Reedy model structure. In this section, we will compare such a Reedy model structure for presheaves with values in Quick’s model category of profinite spaces with the model structures constructed in the previous section. In the second part of this section, we will discuss a new generalized Reedy structure on the category $\Omegaor$ of reduced open trees of a quite different kind. The virtue of this new structure is that it interacts well with filtrations of $\infty$-operads given by bounding the arities of the operations, as we will see in \cref{sec:Boavida-Horel-Robertson}.

\subsection{Comparing \texorpdfstring{\cref{theorem:Main-theorem-constructing-model-structures}}{the main theorem} with the Reedy model structure}

There are two natural candidates for a Reedy-type model structure on $\dProfS$: First, one can consider the model category $\dProfS_R$ obtained by applying \cref{theorem:Main-theorem-constructing-model-structures} to the usual Reedy model structure on $\dSpace$. Secondly, one can use Quick's model structure for profinite $G$-spaces (cf.\ \cref{theorem:Model-Structure-Profinite-G-Spaces-Quick}) and the fact that $\OmegaA$ is a generalized Reedy category in the sense of \cite{BergerMoerdijk2011Extension} to endow $\dProfS$ with a Reedy model structure, as stated in \cref{theorem:True-Reedy-Model-Structure} below (see \cite[\S\S 10.1-10.3]{HeutsMoerdijk2022Trees} for an introduction to (generalized) Reedy categories). The aim of this section is to prove that these model structures agree. In particular, we deduce that the weak equivalences in $\dProfS_R$ are those maps that are levelwise weak equivalences in Quick's model structure for profinite spaces $\s\wh\Set_Q$. This fact will be used to characterize the weak equivalences between profinite $\infty$-operads as the Dwyer--Kan equivalences in \cref{sec:Dwyer-Kan}.

Recall from \cite[\S 1]{BergerMoerdijk2011Extension} that for a dendroidal space $X$ and a tree $T$, the \emph{latching object} of $X$ at $T$ is defined as
\begin{equation*}L_T X = \colim_{T \twoheadrightarrow S} X(S),\end{equation*}
where the colimit is taken over all degeneracies in $\OmegaA$ with source $T$. Dually, the \emph{matching object} of $X$ at the tree $T$ is defined as the limit
\begin{equation*}M_T X = \lim_{S \rightarrowtail T} X(S)\end{equation*}
over all face maps in $\OmegaA$ with target $T$. Note that the latching and matching objects come with natural maps
\begin{equation*}L_T X \to X(T) \to M_T X.\end{equation*}

Let $\eqswhSet{G}_Q$ denote the category of profinite $G$-spaces, endowed with Quick's model structure (see \cref{theorem:Model-Structure-Profinite-G-Spaces-Quick}).

\begin{theorem}\label{theorem:True-Reedy-Model-Structure}
	There exists a simplicial model structure on $\dProfS$ in which a map $X \to Y$ is
	\begin{itemize}[noitemsep]
		\item a cofibration if for every tree $T$, the latching map $X(T) \cup_{L_T X} L_T Y \to Y(T)$ is a normal monomorphism,
		\item a weak equivalence if for every tree $T$, the map $X(T) \to Y(T)$ is a weak equivalence in $\s\wh\Set_Q$, and
		\item a fibration if for every tree $T$, the matching map $X(T) \to M_TX \times_{M_T Y} Y(T)$ is a fibration in $\eqswhSet{\Aut(T)}_Q$.
	\end{itemize}
\end{theorem}

\begin{proof}
	We wish to apply \cite[Thm.\ 1.6]{BergerMoerdijk2011Extension} to the model structures $\eqswhSet{G}_Q$ of \cref{theorem:Model-Structure-Profinite-G-Spaces-Quick}, where $G$ ranges over the groups $\Aut(T)$ for $T \in \OmegaA$. Since we don't know whether Quick's model structure on $\eqswhSet{G}$ is the projective one, the assumptions of that theorem might not be satisfied. However, it is easily checked that the proof still goes through: the only places where the fact that one is working with the projective model structure is used, are in the proofs of Lemmas 5.3 and 5.5 of op.\ cit. One can verify that those proofs go through without a change when using Quick's model structure $\eqswhSet{G}_Q$, hence the result follows.
\end{proof}

\begin{proposition}\label{proposition:Normal-cofibrations-in-profinite-Reedy-model-structure}
	A map in $\dProfS$ is a cofibration in the model structure of \cref{theorem:True-Reedy-Model-Structure} if and only if it is a normal monomorphism.
\end{proposition}

\begin{proof}
	Let $U$ denote the functor that takes a profinite dendroidal space to its underlying dendroidal space. A map $X \to Y$ in $\dProfS$ is a normal monomorphism if and only if the underlying map $UX \to UY$ is a normal monomorphism in $\dSpace$. According to \cite[Lem.\ 12.1]{HeutsMoerdijk2022Trees}, this is the case if and only if for every tree $T$, the latching map
	\begin{equation*}UX(T) \cup_{L_T UX} L_T UY \to UY(T)\end{equation*}
	is injective and $\Aut(T)$ acts freely on the complement of its image. In particular, the result follows if we can show that the underlying set of $X(T) \cup_{L_T X} L_T Y$ agrees with $UX(T) \cup_{L_T UX} L_T UY$. To see that this is the case, note that the relevant colimits are computed by first constructing them in $\Top$ and then applying the functor $L$ that is left adjoint to the inclusion $\Stone \hookrightarrow \Top$ of the category of Stone spaces into all topological spaces. One can show that when these colimits are computed in $\Top$, they are already compact Hausdorff and totally disconnected, hence the left adjoint $L$ does not change them. In particular, the forgetful functor $\Stone \to \Set$ preserves these particular colimits, showing that $X(T) \cup_{L_T X} L_T Y = UX(T) \cup_{L_T UX} L_T UY$.
\end{proof}

\begin{proposition}\label{proposition:Weak-Equivalences-in-profinite-Reedy-model-structure}
	A map $X \to Y$ is a weak equivalence in the model structure $\dProfS_R$ of \cref{theorem:Main-theorem-constructing-model-structures} if and only if for every $T$ in $\OmegaA$, the map $X(T) \to Y(T)$ is a weak equivalence in $\s\wh\Set_Q$.
\end{proposition}

\begin{proof}
	Recall that $X \to Y$ is a weak equivalence in $\dProfS_R$ if and only if
	\begin{equation*}\sHomKan(Y \times E, A) \to \sHomKan(X \times E, A)\end{equation*}
	is a weak equivalence for every lean dendroidal space $A$ that is fibrant in $\dSpace_R$. If $K$ is a lean Kan complex, then the dendroidal space $R_TK$ obtained as the following right Kan extension
	\begin{equation*}\begin{tikzcd}
		* \ar[r,"K"] \ar[d, "T"'] & \s\Set \\
		\OmegaA \ar[ur, dashed, "R_TK"'] &
	\end{tikzcd}\end{equation*}
	is a Reedy fibrant lean dendroidal space. Since $\sHomKan(Z, R_TK) \cong \sHom(Z(T),K)$ for any dendroidal space $Z$, where $\sHom$ denotes the simplicial enrichment of $\s\wh\Set$, it follows that for every tree $T$ and every lean Kan complex $K$, the map
	\begin{equation*}\sHom(Y(T) \times E(T), K) \to \sHom(X(T) \times E(T), K)\end{equation*}
	is a weak equivalence. In particular, $X(T) \times E(T) \to Y(T) \times E(T)$ is a weak equivalence in $\s\wh\Set_Q$. Since $E(T) \to *$ is a trivial fibration by construction, it follows that $X(T) \to Y(T)$ is a weak equivalence in $\s\wh\Set_Q$.
	
	For the converse, it suffices to show that any lean dendroidal space $A$ that is fibrant in $\dSpace_R$ is also fibrant in the model structure of \cref{theorem:True-Reedy-Model-Structure}. Namely, since the latter model structure is simplicial, this would imply that
	\begin{equation*}\sHomKan(Y \times E, A) \to \sHomKan(X \times E, A)\end{equation*}
	is a weak equivalence for any pointwise equivalence $X \to Y$. By definition, $A$ is fibrant in the model structure of \cref{theorem:True-Reedy-Model-Structure} if for every tree $T$, the map
	\begin{equation}\label{equation:Matching-map}
		A(T) \to M_T A
	\end{equation}
	is a fibration in the model structure $\eqswhSet{G}_Q$ of \cref{theorem:Main-theorem-constructing-model-structures}, where $G$ denotes $\Aut(T)$. Since $A$ is lean and fibrant in $\dSpace_R$, the map \cref{equation:Matching-map} is a fibration in the projective model structure on $\eqsSet{G}$ between lean $G$-Kan complexes, hence this is indeed the case by \cref{corollary:Fibrations-between-lean-profinite-G-spaces}.
\end{proof}

\begin{corollary}
	The model structure of \cref{theorem:True-Reedy-Model-Structure} agrees with $\dSpace_R$.
\end{corollary}

\begin{proof}
	This follows since the cofibrations and the weak equivalences of both model structures agree by \cref{proposition:Normal-cofibrations-in-profinite-Reedy-model-structure,proposition:Weak-Equivalences-in-profinite-Reedy-model-structure}.
\end{proof}

\subsection{The outer Reedy structure on \texorpdfstring{$\Omegaor$}{Ω\_or}}\label{ssec:Outer-Reedy-structure}

In \cref{sec:Boavida-Horel-Robertson}, we will make extensive use of a certain filtration
\begin{equation}\label{eq:filtration-Omegaor-by-arity}
	\{\eta\} = \Omegaorf{1} \subset \Omegaorfp{1} \subset \Omegaorf{2} \subset \Omegaorfp{2} \subset \cdots \subset \cup_n \Omegaorf{n} = \Omegaor
\end{equation}
of $\Omegaor$. This filtration is inspired by a filtration introduced by Göppl \cite{Goppl2019SpectralSequenceSpaces}, which he used to analyse the connectivity of the derived mapping spaces between certain operads. The category $\Omegaorf{n}$ is defined as the full subcategory of $\Omegaor$ spanned by those trees where each vertex has at most $n$ ingoing edges, and $\Omegaorfp{n}$ denotes the full subcategory of $\Omegaor$ that has one extra object, namely the $(n+1)$-corolla $C_{n+1}$. An issue with this filtration is that the restriction and right Kan extension functors corresponding to these inclusions are not always right Quillen with respect to the Reedy model structures corresponding to the usual Reedy structure on $\Omegaor$. We will remedy this by considering an alternative Reedy structure on $\Omegaor$, called the \emph{outer Reedy structure}, which behaves better with respect to this filtration (cf.\ \cref{remark:Reason-for-using-outer-Reedy}). Throughout this section, we will refer to the usual Reedy structure on $\Omegaor$ as the \emph{standard Reedy structure}. Recall the definitions of outer and inner face maps from \cite[\S 3.3]{HeutsMoerdijk2022Trees}. We will simply call a composition of outer face maps an outer face map, and similarly for (compositions of) inner face maps.

\begin{definition}
	In the \emph{outer Reedy structure} on $\Omegaor$, the positive morphisms are defined to be the outer face maps, denoted $\oface$, while the negative morphisms are defined to be the inner face maps, denoted $\iface$. The degree function for this Reedy structure is given by sending a tree $T$ to its \emph{weight}, which is defined by
	\begin{equation}\label{equation:Definition-of-weight-of-a-tree}
		w(T) = 2 \cdot \# l(T) - \# v(T).
	\end{equation}
	Here $l(T)$ denotes the set of leaves of $T$ and $v(T)$ its set of vertices.
\end{definition}

We call this the \emph{outer} Reedy structure since matching objects with respect to this Reedy structure are computed in terms of outer face maps, as opposed to all face maps. Note that the outer Reedy structure is indeed a (dualizable, generalized) Reedy structure on $\Omegaor$ in the sense of \cite[Def.\ 1.1]{BergerMoerdijk2011Extension} and \cite[Def.\ 10.1]{HeutsMoerdijk2022Trees}: First, it is a straightforward check that $w$ defines a degree function. Moreover, it follows from \cite[Prop.\ 3.10]{HeutsMoerdijk2022Trees} that any map in $\Omegaor$ factors uniquely up to unique isomorphism as an inner face map followed by an outer face map. Finally, item (d) of \cite[Def.\ 10.1]{HeutsMoerdijk2022Trees} follows since any automorphism in $\Omegaor$ is completely determined by what it does on leaves.

The outer Reedy structure on $\Omegaor$ restricts to Reedy structures on $\Omegaorf{n}$ and $\Omegaorfp{n}$, which we will also call outer Reedy structures. The categories of $\s\wh\Set$-valued presheaves on $\Omegaorf{n}$ and $\Omegaorfp{n}$ are denoted $\rodProfSf{n}$ and $\rodProfSfp{n}$, respectively. Given a presheaf $X$ on $\Omegaor$, $\Omegaorf{n}$ or $\Omegaorfp{n}$, we write
\begin{equation*}
	\Lo_T X = \colim_{T \iface S} X(S) \quad \text{and} \quad \Mo_T X = \lim_{S \oset[-0.6ex]{\circ}{\longrightarrow} T} X(S)
\end{equation*}
for the latching and matching objects at $T$ with respect to the outer Reedy structure.

The outer Reedy structure on $\Omegaor$ has the property that for any $n \geq 0$, the number of (isomorphism classes of) trees with weight at most $n$ is finite. In particular, one obtains the following result by the same proof as \cref{lemma:Detecting-lean-using-matching-maps}.

\begin{lemma}\label{lemma:Detecting-lean-using-outer-matching-maps}
	Let $X$ be a dendroidal space. Then $X$ is lean if and only if 
	\begin{enumerate}[(1)]
		\item for every tree $T$, the simplicial set $X(T)$ is lean, and
		\item there exists an $n$ such that for every tree $T$ of weight $w(T) > n$, the matching map $X \to \Mo_T X$ is an isomorphism. \qed
	\end{enumerate}
\end{lemma}

Exactly as in \cref{theorem:True-Reedy-Model-Structure}, one obtains Reedy model structures on $\rod\ProfS$, $\rodProfSf{n}$ and $\rodProfSfp{n}$ with respect to the outer Reedy structures. These will be denoted $\rod\ProfS_{oR}$, $\rodProfSf{n}_{oR}$ and $\rodProfSfp{n}_{oR}$, respectively. A fibrant object in one of these model categories will be called \emph{outer fibrant}. Note that the outer Reedy model structure has more fibrations and fewer cofibrations than the standard Reedy model structure. The dual of the proof of \cite[Thm.\ 10.14]{HeutsMoerdijk2022Trees} shows that these model structures are fibrantly generated, where the generating (trivial) fibrations are maps between lean outer fibrant objects.

We wish to localize these to obtain model structures analogous to $\rod\ProfS_{RS}$. For this we will use the following lemma. By a lean outer fibrant reduced open dendroidal Segal space, we mean a lean outer fibrant object in $\rod\Space$ satisfying the Segal condition. Note that this slightly conflicts with the terminology of \cref{ssec:Model-structures-preliminaries}, where we defined reduced open dendroidal Segal spaces as the fibrant objects of $\rod\Space_{RS}$, implicitly assuming they are Reedy fibrant.

\begin{lemma}\label{lemma:Lean-outer-Reedy-Segal-spaces-vs-standard-Reedy}
	Let $X \in \rod\Space$ be a lean outer fibrant reduced open dendroidal Segal space. Then there exists a lean standard Reedy fibrant object $Y$ that is pointwise equivalent to $X$.
\end{lemma}

\begin{proof}
	The same argument as in the standard Reedy case (cf.\ the first line of the proof of \cref{corollary:Characterizing-lean-fibrant-reduced}) shows that $X$ is $\pi$-finite. Let $X \wearrow \wt X$ be a standard Reedy fibrant replacement of $X$ in $\rod\Space$. Then $\wt X$ is $\pi$-finite and satisfies the Segal condition, hence by \cref{corollary:Characterizing-lean-fibrant-reduced} it follows that there exists a lean standard Reedy fibrant reduced open dendroidal Segal space $Y$ together with a Reedy equivalence $\wt X \wearrow Y$.
\end{proof}

Using the work of Christensen and Isaksen on right Bousfield localizations \cite[Thm.\ 2.6]{IsaksenChristensen2004ProSpectra}, we can now deduce the following.

\begin{proposition}\label{proposition:Existence-outer-Reedy-Segal-model-structure}
	There exists a model structure $\rod\ProfS_{oRS}$ that is a left Bousfield localization of the outer Reedy model structure on $\rod\ProfS$ and whose weak equivalences agree with those of $\rod\ProfS_{RS}$. Any fibrant object in this model category is a cofiltered limit of lean outer fibrant reduced open dendroidal Segal spaces.
\end{proposition}

\begin{proof}
	Let $T$ denote the collection of lean outer fibrant reduced open dendroidal Segal spaces. By \cref{lemma:Lean-outer-Reedy-Segal-spaces-vs-standard-Reedy}, the $T$-local equivalences in $\rod\ProfS$ are precisely the weak equivalences of $\rod\ProfS_{RS}$. Note that since the category of lean reduced open dendroidal spaces is essentially small, we may assume $T$ to be a set. By the dual of \cite[Thm.\ 2.6]{IsaksenChristensen2004ProSpectra}, the left Bousfield localization $\rod\ProfS_{oRS}$ of $\rod\ProfS_{oR}$ in which the weak equivalences are the $T$-local ones exists.
	
	For the statement about fibrant objects in this model category, note that one can construct fibrant replacements in $\rod\ProfS_{oRS}$ using the cosmall object argument by the dual of \cite[Lem.\ 2.5]{IsaksenChristensen2004ProSpectra}. Using the dual of the fat small object argument of \cite{MakkaiRosickyea2014FatSmallObject} instead, it follows that any fibrant object is a cofiltered limit of objects $A$ with the property that $A \to *$ is a finite composition of pullbacks of maps such that each of these maps is either of the form $B^{\Delta[n]} \to B^{\partial \Delta[n]}$, where $B$ is an object of $T$, or is a generating fibration in $\rod\ProfS_{oR}$. Since such maps are fibrations between lean outer fibrant reduced open dendroidal Segal spaces, it follows that $A$ itself is a lean outer fibrant reduced open dendroidal Segal space. Given a fibrant object $X$ in $\rod\ProfS_{oRS}$, this argument gives us a trivial cofibration $X \trivcofarrow \wt X$ where $\wt X$ is an inverse limit of lean outer fibrant reduced open dendroidal Segal spaces. Since both $X$ and $\wt X$ are fibrant, $X$ is a retract of $\wt X$ and must itself be an inverse limit of lean outer fibrant reduced open dendroidal Segal spaces.
\end{proof}

\begin{remark}
	Note that one could also have constructed the model structure $\gd\ProfS_{RSC}$ from \cref{theorem:Main-theorem-constructing-model-structures} using the method given here. However, it would not be clear from this construction that one obtains a fibrantly generated model category, and results such as \cref{theorem:Profinite-Segal-spaces-vs-profinite-operads} do not seem to follow as easily from this construction.
\end{remark}

Our proofs in \cref{sec:Boavida-Horel-Robertson} also require analogous model structures on $\rodProfSf{n}$ and $\rodProfSfp{n}$. We will say that an object $X$ in one of these categories satisfies the \emph{Segal condition} if for any tree $R = S \circ_e T$ where $e$ is an internal edge of $R$, the map
\begin{equation*}
	X(R) \to X(S) \times_{X(e)}^h X(T)
\end{equation*}
is a weak equivalence. Write $T_n$ and $T_n^+$ for the collection of lean outer fibrant objects in $\rodProfSf{n}$ and $\rodProfSfp{n}$, respectively. The proof of \cref{proposition:Existence-outer-Reedy-Segal-model-structure} now yields the following.

\begin{proposition}\label{proposition:Existence-outer-Reedy-Segal-model-structure-restricted}
	There exist left Bousfield localizations $\rodProfSf{n}_{oRS}$ and $\rodProfSfp{n}_{oRS}$ of the outer Reedy model structures on $\rodProfSf{n}$ and $\rodProfSfp{n}$ in which the weak equivalences are the $T_n$ and $T_n^+$-local ones, respectively. Moreover, any fibrant object in $\rodProfSf{n}_{oRS}$ is a cofiltered limit of objects in $T_n$ and any fibrant object in $\rodProfSfp{n}_{oRS}$ is a cofiltered limit of objects in $T_n^+$. \qed
\end{proposition}

Let us conclude this section by describing a strategy to detect fibrant objects in these model structures. This will be used in \cref{sec:Boavida-Horel-Robertson}, where we show that an object is fibrant in $\rod\ProfS_{oRS}$ if and only if it is outer fibrant and satisfies the Segal condition. The following standard fact on left Bousfield localizations will be used several times. Note that this result does not follow automatically from the way these Bousfield localizations are constructed in \cite[\S 2]{IsaksenChristensen2004ProSpectra}.

\begin{lemma}\label{lemma:Detecting-fibrant-objects-in-left-Bousfield-localization}
	Let $L\mathcal{M}$ be a left Bousfield localization of a model category $\mathcal{M}$ and let $X$ be fibrant in $L\mathcal{M}$. If $Y \wearrow X$ is a weak equivalence in $\mathcal{M}$ and $Y$ is fibrant in $\mathcal{M}$, then $Y$ is fibrant in $L\mathcal{M}$.
\end{lemma}

\begin{proof}
	Factor $Y \wearrow X$ in $\mathcal{M}$ as a trivial cofibration $Y \trivcofarrow \wt Y$ followed by a trivial fibration $\wt Y \trivfibarrow X$. Since $Y \trivcofarrow \wt Y$ is a trivial cofibration between fibrant objects in $\mathcal{M}$, it admits a retract. Moreover, as the trivial fibrations of $\mathcal{M}$ and $L\mathcal{M}$ agree, it follows that $\wt Y$ is fibrant in $L\mathcal{M}$. Since $Y$ is a retract of a fibrant object in $L\mathcal{M}$, it must itself be fibrant in $L\mathcal{M}$.
\end{proof}

One can deduce from this that the full subcategory of fibrant objects in $\rod\ProfS_{oRS}$ is closed under those cofiltered limits which exist in the full subcategory of fibrant objects in $\rod\ProfS_{oR}$.

\begin{lemma}\label{lemma:Detecting-fibrant-objects}
	Let $X$ be an outer fibrant object in $\rod\ProfS$. Then the following are equivalent:
	\begin{enumerate}[(i),noitemsep]
		\item\label{item1:lemma:Detecting-fibrant-objects} $X$ is fibrant in $\rod\ProfS_{oRS}$.
		\item\label{item2:lemma:Detecting-fibrant-objects} $X$ is a cofiltered limit of objects that are fibrant in $\rod\ProfS_{oRS}$.
		\item\label{item3:lemma:Detecting-fibrant-objects} $X$ is a cofiltered limit of lean outer fibrant objects satisfying the Segal condition.
	\end{enumerate}
	The analogous statement holds for $\rodProfSf{n}$ and $\rodProfSfp{n}$.
\end{lemma}

\begin{proof}
	Let $\mathcal{N}_{oRS}$ denote one of the model categories $\rod\ProfS_{oRS}$, $\rodProfSf{n}_{oRS}$ and $\rodProfSfp{n}_{oRS}$. It follows from \cref{proposition:Existence-outer-Reedy-Segal-model-structure,proposition:Existence-outer-Reedy-Segal-model-structure-restricted} that \ref{item1:lemma:Detecting-fibrant-objects} implies \ref{item3:lemma:Detecting-fibrant-objects}. Since lean outer fibrant objects satisfying the Segal condition are by construction fibrant in $\mathcal{N}_{oRS}$, item \ref{item2:lemma:Detecting-fibrant-objects} follows from \ref{item3:lemma:Detecting-fibrant-objects}. Finally, to deduce \ref{item1:lemma:Detecting-fibrant-objects} from \ref{item2:lemma:Detecting-fibrant-objects}, suppose that $X = \lim_{i \in I} X_i$ is an outer fibrant object that is a cofiltered limit of objects that are fibrant in $\mathcal{N}_{oRS}$. Then let $\{\wt X_i\}_{i \in I}$ be a fibrant replacement of $\{X_i\}_{i \in I}$ in the injective model structure on $\Fun(I^{op},\mathcal{N}_{oRS})$. Since $\lim_i \colon \Fun(I^{op}, \mathcal{N}_{oRS}) \to \mathcal{N}_{oRS}$ is right Quillen with respect to the injective model structure, we see that $\lim_i \wt X_i$ is fibrant in $\mathcal{N}_{oRS}$. The maps $X_i \to \wt X_i$ are weak equivalences between fibrant objects in $\mathcal{N}_{oRS}$, hence by \cite[Thm.\ 3.2.13]{Hirschhorn2003Model} they are Reedy equivalences. Since these are stable under cofiltered limits by \cref{proposition:Weak-equivalences-stable-under-cofiltered-limits}, we see that $X = \lim_i X_i \to \lim_i \wt X_i$ is a Reedy equivalence. By \cref{lemma:Detecting-fibrant-objects-in-left-Bousfield-localization}, we conclude that $X$ is fibrant in $\mathcal{N}_{oRS}$.
\end{proof}
	
\end{document}

%% file: Sections/DwyerKanEquivalences.tex
	\section{Weak equivalences of profinite \texorpdfstring{$\infty$}{infinity}-operads}\label{sec:Dwyer-Kan}
	
	The weak equivalences between ordinary $\infty$-operads can be characterized as those maps that are surjective up to homotopy equivalence and that induce weak homotopy equivalences on the spaces of operations. In this section, we, extend this characterization to the setting of profinite $\infty$-operads and profinite complete dendroidal Segal spaces.
	
	We start by proving that weak equivalences are detected by the forgetful functor to $\dSpace$. The main ingredient is the following theorem, which follows from Theorem E.3.1.6 of \cite{Lurie2016SAG} or the second part of the proof of \cite[Prop.\ 3.9]{BoavidaHorelRobertson2019GenusZero}.
	
	\begin{theorem}\label{theorem:Profinite-spaces-equivalence-detected-underlying}
		Let $X \to Y$ be a map between fibrant objects in the Quick model structure on $\s\wh\Set$. Then $X \to Y$ is a weak equivalence if and only if the underlying map $UX \to UY$ is a weak equivalence in the Kan--Quillen model structure on $\s\Set$.
	\end{theorem}

	This implies the following:
	
	\begin{proposition}\label{proposition:Weak-equivalence-profinite-dendroidal-spaces-underlying}
		Let $\dSpace$ be endowed with one of the model structures of \cref{ssec:Model-structures-preliminaries} and let $\dProfS$ be endowed with the corresponding model structure of \cref{theorem:Main-theorem-constructing-model-structures}. Then a map $X \to Y$ between fibrant objects in $\dProfS$ is a weak equivalence if and only if the underlying map $UX \to UY$ is a weak equivalence in $\dSpace$.
	\end{proposition}

	\begin{proof}
		One direction follows since $U \colon \dProfS \to \dSpace$ is right Quillen. For the other direction, note that the model structure on $\dProfS$ is a left Bousfield localization of $\dProfS_R$ (since their cofibrations agree). By \cref{proposition:Weak-Equivalences-in-profinite-Reedy-model-structure} and \cite[Thm.\ 3.2.13]{Hirschhorn2003Model}, we see that a map $X \to Y$ between fibrant objects in $\dProfS$ is a weak equivalence if and only if it is a pointwise weak equivalence in Quick's model structure $\s\wh\Set_Q$. By \cref{theorem:Profinite-spaces-equivalence-detected-underlying}, this is the case if and only if $UX(T) \to UY(T)$ is a weak homotopy equivalence for every tree $T$, which follows if $UX \to UY$ is a weak equivalence.
	\end{proof}

	\begin{corollary}\label{corollary:Weak-equivalences-profinite-infty-operads-underlying}
		Let $\OmegaA = \Omega, \Omegao, \Omegacl$. Then a map $X \to Y$ is a weak equivalence between fibrant objects in the operadic model structure on $\dwhSet$ if and only if the underlying map $UX \to UY$ is a weak equivalence in the operadic model structure on $\dSet$.
	\end{corollary}

	\begin{proof}
		This follows from \cref{proposition:Weak-equivalence-profinite-dendroidal-spaces-underlying}, using the Quillen equivalences \cref{equation:Quillen-equivalence-dSpace-dSet-Sing}-\cref{equation:Quillen-equivalence-dSpace-dSet-0} and \cref{theorem:Profinite-Segal-spaces-vs-profinite-operads}.
	\end{proof}

	One can mimic the usual definition of the spaces of operations given in \cref{ssec:Spaces-of-operations} in the profinite setting. Given a profinite dendroidal Segal space $X$, its \emph{colours} are defined as the maps $* \to X$; that is, the points of $X(\eta)_0$ when viewed as a Stone space. Given a tuple $c_1,\ldots,c_n, d \in X(\eta)_0$ of colours, we define the \emph{(profinite) space of operations} $X(c_1,\ldots,c_n;d)$ as the pullback
	\begin{equation*}\begin{tikzcd}
		X(c_1,\ldots,c_n;d) \ar[r] \ar[d] \ar[dr, phantom, very near start, "\lrcorner"] & X(C_n) \ar[d] \\
		\Delta[0] \ar[r,"(c_1 {,} \ldots {,} c_n {,} d)"] & X(\eta)^{n+1}.
	\end{tikzcd}\end{equation*}
	In the closed case, one should replace $\eta$ and $C_n$ with their closures $\ol{\eta}$ and $\ol{C_n}$. Note that the Reedy fibrancy of $X$ implies that the right-hand vertical map is a fibration in $\s\wh\Set_Q$, hence $X(c_1,\ldots,c_n;d)$ is fibrant in Quick's model structure. Given two colours $c,d \in X(\eta)_0$, we say that $c$ and $d$ are \emph{homotopy equivalent} if there exists a map $f \colon J \to X$ such that $f(0)= c$ and $f(1) = d$. If $X$ is a complete profinite dendroidal Segal space, then by completeness this is equivalent to the existence of a $1$-simplex $x \in X(\eta)_1$ such that $d_1x = c$ and $d_0 x = d$. This leads to the following definition.
	
	\begin{definition}
		A map $f \colon X \to Y$ of complete profinite dendroidal Segal spaces is \emph{fully faithful} if for any tuple $c_1,\ldots,c_n,d \in X(\eta)_0$ of colours, the map $X(c_1,\ldots,c_n;d) \to Y(fc_1,\ldots,fc_n;fd)$ is a weak equivalence in $\s\wh\Set_Q$, and it is \emph{essentially surjective} if for every colour $d \in Y(\eta)_0$, there exists a $c \in X(\eta)_0$ such that $fc$ is homotopy equivalent to $d$. A map that is both fully faithful and essentially surjective is called a \emph{Dwyer--Kan equivalence}.
	\end{definition}

	As before, one should replace $\eta$ with its closure $\ol \eta$ in this definition in the closed case.
	
	\begin{theorem}\label{theorem:Dwyer-Kan-equivalence-Segal-case}
		A map $X \to Y$ of complete profinite dendroidal Segal spaces is a weak equivalence in $\dProfS_{RSC}$ if and only if it is a Dwyer--Kan equivalence.
	\end{theorem}

	\begin{proof}
		It is clear that two colours $c,d \in Y(\eta)_0$ are homotopy equivalent if and only if they are homotopy equivalent as colours of $UY$. In particular, $X \to Y$ is essentially surjective if and only if $UX \to UY$ is. Moreover, since $U \colon \s\wh\Set \to \s\Set$ preserves pullbacks, we see that $U(X(c_1,\ldots,c_n;d)) \cong (UX)(c_1,\ldots,c_n;d)$ and similarly for $Y$. Combined with \cref{theorem:Profinite-spaces-equivalence-detected-underlying}, we see that $X \to Y$ is a Dwyer--Kan equivalence if and only if the underlying map $UX \to UY$ is a Dwyer--Kan equivalence of complete dendroidal Segal spaces. By \cite[Cor.\ 12.34]{HeutsMoerdijk2022Trees}, this is the case if and only if $UX \to UY$ is a weak equivalence, which by \cref{proposition:Weak-equivalence-profinite-dendroidal-spaces-underlying} is the case if and only if $X \to Y$ is a weak equivalence.
	\end{proof}

	We aim to also define Dwyer--Kan equivalences for maps between fibrant objects in $\dwhSet$. This is slightly more involved: if we mimic the usual definition for the spaces of operations, then we obtain ordinary simplicial sets. However, if $Z$ is a degreewise finite dendroidal set and $X$ is a dendroidal profinite set, then the $\sHom(Z,X)$ can be canonically upgraded to a simplicial profinite set. It turns out that if $X$ is an $\infty$-operad, then this simplicial profinite set is a \emph{profinite $\infty$-category}, i.e.\ a fibrant object in the profinite Joyal model structure $\s\wh\Set_J$ introduced in Example 5.7 of \cite{BlomMoerdijk2020SimplicialProPublished}. To see this, first observe that for a lean dendroidal set $A$, the simplicial set $\sHom(Z,A)$ is lean by \cite[Lem.\ 2.9]{BlomMoerdijk2022ProfinitePartI}. This implies $\sHom(-,-)$ extends uniquely to a functor
	\begin{equation*}
		\wh\sHom(-,-) \colon \gd\Fin\Set^{op} \times \dwhSet \to \s\wh\Set; \quad (Z,\{X_i\}_{i \in I}) \mapsto \lim_i \sHom(Z,X_i)
	\end{equation*}
	that preserves cofiltered limits in the second variable. Moreover, the underlying simplicial set $U(\wh\sHom(Z,X))$ is canonically isomorphic to $\sHom(Z,X)$. The compatibility of $\wh\sHom(-,-)$ with the profinite Joyal model structure $\s\wh\Set_J$ now takes the following form. 
	
	\begin{lemma}
		Let $i \colon Z \cofarrow W$ be a normal monomorphism between normal degreewise finite dendroidal sets and let $p \colon Y \fibarrow X$ be a fibration in the operadic model structure on $\dwhSet$. Then
		\begin{equation}\label{equation:Pullback-power-profinite-hom}
			\wh\sHom(W,Y) \to \wh\sHom(Z,Y) \times_{\wh\sHom(Z,X)} \wh\sHom(W,X)
		\end{equation}
		is a fibration in the profinite Joyal model structure $\s\wh\Set_J$ of \cite[Ex.\ 5.7]{BlomMoerdijk2020SimplicialProPublished}, which is trivial if either $i$ or $p$ is.
	\end{lemma}

	\begin{proof}
		It suffices to check this when $Y \fibarrow X$ is a generating (trivial) fibration, i.e.\ an operadic fibration between lean dendroidal sets. In this case it follows from \cite[Prop.\ 9.28]{HeutsMoerdijk2022Trees} that the map \cref{equation:Pullback-power-profinite-hom} is a (trivial) fibration between lean quasi-categories in the Joyal model structure on $\s\Set$, hence that it is a generating (trivial) fibration of $\s\wh\Set_J$.
	\end{proof}

	One now defines the \emph{colours} of a profinite $\infty$-operad $X$ as the points of the Stone space $X(\eta)$. 
	Given a tuple $c_1,\ldots,c_n, d \in X(\eta)$ of colours, we define the \emph{(profinite) space of operations} $X(c_1,\ldots,c_n;d)$ as the pullback
	\begin{equation*}\begin{tikzcd}
		X(c_1,\ldots,c_n;d) \ar[r] \ar[d] \ar[dr, phantom, very near start, "\lrcorner"] & \wh\sHom(\Omega[C_n],X) \ar[d] \\
		\Delta[0] \ar[r,"(c_1 {,} \ldots {,} c_n {,} d)"] & \wh\sHom(\Omega[\eta],X)^{n+1}.
	\end{tikzcd}\end{equation*}
	Again, $\eta$ and $C_n$ need to be replaced with their closures in the case of closed dendroidal sets. 
	Note that the right vertical map is a fibration in the profinite Joyal model structure, hence $X(c_1,\ldots,c_n;d)$ is a profinite $\infty$-category.\footnote{In fact, one can show that $X(c_1,\ldots,c_n;d)$ is fibrant in the Quick model structure $\s\wh\Set_Q$. However, since the proof is somewhat technical and we don't need it here, we omit it.} 
	Given two colours $c,d \in X(\eta)$, we say that $c$ and $d$ are \emph{homotopy equivalent} if there exists a map $f \colon J \to X$ such that $f(0) = c$ and $f(1) = d$.
	
	\begin{definition}
		A map $f \colon X \to Y$ of profinite $\infty$-operads is \emph{fully faithful} if for any tuple $c_1,\ldots,c_n,d \in X(\eta)_0$ of colours, the map $X(c_1,\ldots,c_n;d) \to Y(fc_1,\ldots,fc_n;fd)$ is a weak equivalence in the profinite Joyal model structure $\s\wh\Set_J$. It is \emph{essentially surjective} if for any $c \in Y(\eta)$, there exists a $d \in X(\eta)$ such that $c$ is homotopy equivalent to $f(d)$. A map that is both fully faithful and essentially surjective is called a \emph{Dwyer--Kan equivalence}.
	\end{definition}
	
	One readily deduces the following analogue of \cref{theorem:Dwyer-Kan-equivalence-Segal-case}.
	
	\begin{theorem}\label{theorem:Dwyer-Kan-equivalence-profinite-infty-operads}
		A map $X \to Y$ of profinite $\infty$-operads is a weak equivalence if and only if it is a Dwyer--Kan equivalence.
	\end{theorem}

	\begin{proof}
		It is clear that if $X \to Y$ is a weak equivalence, then it is a Dwyer--Kan equivalence. Conversely, if it is a Dwyer--Kan equivalence, then the underlying map $UX \to UY$ is also a Dwyer--Kan equivalence since $U$ commutes with the pullbacks used to define the spaces of operations. By \cite[Thm.\ 3.11]{CisinskiMoerdijk2013Segal}, $UX \to UY$ is an equivalence of $\infty$-operads, hence by \cref{corollary:Weak-equivalences-profinite-infty-operads-underlying} we conclude that $X \to Y$ is a weak equivalence.
	\end{proof}

%% file: Sections/UnderlyingInftyCategory.tex
 	\section{The underlying \texorpdfstring{$\infty$}{infinity}-category}\label{sec:Underlying-infty}
 	
 	The aim of this section is to identify the underlying $\infty$-categories of $\dProfS_{RSC}$ and $\dwhSet$. Recall the meaning of the subscript ``$RSC$'' from \cref{notation:Subscripts-dProfS}. In light of the Quillen equivalence between $\dProfS_{RSC}$ and the operadic model structure on $\dwhSet$ in the case that $\OmegaA = \Omega, \Omegao, \Omegacl$, we also obtain the description of the underlying $\infty$-category of $\dwhSet$ announced in \cite[Thm.\ B]{BlomMoerdijk2022ProfinitePartI}.
 	
 	A possible approach to proving this result would be to try and mimic the arguments of \cite[Appendix]{BlomMoerdijk2020SimplicialProPublished}. However, this requires some care since not all objects of $\dProfS_{RSC}$ are cofibrant. We will circumvent this by constructing a Quillen equivalence between $\dProfS_{RSC}$ and an alternative model category to which \cite[Thm.\ A.7]{BlomMoerdijk2020SimplicialProPublished} applies directly.
 	
 	\subsection{Slicing over a normalization of the terminal object}
 	
 	Let $E \trivfibarrow *$ in $\dSpace$ be a degreewise finite normalization in the sense of \cref{lemma:Degreewise-finite-normalization-point}. Viewing it as an object of $\dProfS$, we can define the slice model category $\dProfS/E$ with respect to any of the model structures of \cref{theorem:Main-theorem-constructing-model-structures}. The map $E \trivfibarrow *$ is a trivial fibration in $\dProfS$ by the same argument as \cite[Lem.\ 5.3]{BlomMoerdijk2022ProfinitePartI}, so we see that
 	\begin{equation*}U : \dProfS/E \rightleftarrows \dProfS : - \times E.\end{equation*}
 	is a Quillen equivalence. Note that we similarly have a Quillen equivalence
 	\begin{equation*}U : \dSpace/E \rightleftarrows \dSpace : - \times E,\end{equation*}
 	where $\dSpace$ is endowed with any of the model structures of \cref{ssec:Model-structures-preliminaries}.
 	
 	In order to apply \cite[Thm.\ A.7]{BlomMoerdijk2020SimplicialProPublished}, we first need to identify $\dProfS/E$ with a pro-category.
 	
 	\begin{definition}
 		Let $\Cin$ be a category and $B \in \Pro(\Cin)$ a pro-object. Define the category $\CinB$ as the full subcategory of $\Pro(\Cin)/B$ consisting of those maps $X \to B$ for which there exists a pullback square of the form
 		\begin{equation*}\begin{tikzcd}
 			X \ar[r] \ar[d] \ar[dr, phantom, near start, "\lrcorner"] & X' \ar[d] \\
 			B \ar[r] & B',
 		\end{tikzcd}\end{equation*}
 		where $B'$ and $X'$ lie in the image of the inclusion $\Cin \hookrightarrow \Pro(\Cin)$.
 	\end{definition}
 
	\begin{remark}\label{remark:Relatively-lean-for-explicit-presentation}
		Let $B = \{B_i\}_{i \in I}$ in $\Pro(\Cin)$ and denote the projection maps by $p_i \colon B \to B_i$. Since any map from $B$ to an object $B'$ of $\Cin$ factors through some $B_i$, we see that $X \to B$ is an object of $\CinB$ if and only if there exists an $i \in I$ such that $X \cong p_i^*X'$ for some $X'$ in $\Cin/B_i$.
	\end{remark}
 	
 	\begin{lemma}
 		Suppose $\Cin$ is a small category that admits finite limits. Then the inclusion $\CinB \hookrightarrow \Pro(\Cin)/B$ induces an equivalence of categories
 		\begin{equation*}\Pro(\CinB) \wearrow \Pro(\Cin)/B\end{equation*}
 	\end{lemma}
 
	\begin{proof}
		We apply \cite[Lem.\ 3.2]{BlomMoerdijk2022ProfinitePartI} to the inclusion $\CinB \hookrightarrow \Pro(\Cin)/B$. To see that all objects of $\CinB$ are cocompact in $\Pro(\Cin)/B$, let $p \colon B \to B'$ be any map from $B$ to an object $B'$ of $\Cin$. Pullback and postcomposition with $p$ give an adjunction
		\begin{equation*}p_! : \Pro(\Cin)/B \rightleftarrows \Pro(\Cin)/B': p^*.\end{equation*}
		Since the left adjoint $p_!$ preserves cofiltered limits, it follows that $p^*$ must preserve cocompact objects. The claim now follows from the equivalence $\Pro(\Cin)/B' \simeq \Pro(\Cin/B')$, which we leave to the reader to verify.
		
		To see that any element of $\Pro(\Cin)/B$ can be written as a cofiltered limit of objects in $\CinB$, let $Y \to B$ be an object of $\Pro(\Cin)/B$ and assume without loss of generality that it has a level representation $\{Y_i \to B_i\}_{i \in I}$ by maps between objects of $\Cin$ (cf.\ Corollary 3.2 of \cite[Appendix]{ArtinMazur1986Etale}). Let $p_i \colon B \to B_i$ denote the structure maps of the diagram $\{B_i\}_{i \in I}$. Then $Y \cong \lim_i p_i^*Y_i$, hence $Y$ is an inverse limit of objects in $\CinB$.
	\end{proof}

	\begin{example}\label{example:Writing-slice-over-E-as-pro-category}
		Applying this to $\dProfS \simeq \Pro(\Leandspace)$ and the object $E = \{\cosk_n E\}_{n \in \bbN}$ constructed in \cref{lemma:Degreewise-finite-normalization-point}, we find that \begin{equation*}\dProfS/E \simeq \Pro(\Leandspace_E).\end{equation*}
		If we write $p_n$ for the projection $E \to \cosk_n E$, then by \cref{remark:Relatively-lean-for-explicit-presentation}, the category $\Leandspace_E$ consists precisely of objects of the form $p_n^* X$, where $X$ is a lean dendroidal space with a map to $\cosk_n E$.
	\end{example}
 	
 	\subsection{Identifying the underlying \texorpdfstring{$\infty$}{infinity}-category}
 	
 	Suppose throughout this section that $\dSpace$ is endowed with any of the model structures of \cref{ssec:Model-structures-preliminaries} and $\dSpace/E$ with the slice model structure. Any $X$ in $\Leandspace_E$ is degreewise finite, hence we may view $\Leandspace_E$ as a full subcategory of $\dSpace/E$. Since any dendroidal space admitting a map to $E$ must be normal, we see that all objects in $\dSpace/E$ are cofibrant. Moreover, the simplicial enrichment $\sHomKan$ on $\dSpace$ induces a simplicial enrichment $\sHom_{/E}$ on $\dSpace/E$, making it into a simplicial model category. In particular, by \cite[Ex.\ 5.3]{BlomMoerdijk2020SimplicialProPublished} $\Leandspace_E$ inherits the structure of a fibration test category in the sense of \cite[Def.\ 5.1]{BlomMoerdijk2020SimplicialProPublished}. Applying Theorem 5.2 of op.\ cit.\ now yields the following result.
 	
 	\begin{proposition}\label{proposition:Completed-model-structure-on-pro-slice}
 		There exists a simplicial model structure on $\Pro(\Leandspace_E)$ in which a map $X \to Y$ is a weak equivalence if and only if
 		\begin{equation*}\sHom_{/E}(Y,A) \to \sHom_{/E}(X,A)\end{equation*}
 		is a weak equivalence for every $A$ in $\Leandspace_E$ that is fibrant in $\dSpace$. It admits sets of generating fibrations and generating trivial fibrations given by
 		\begin{equation*}\mathcal{P} = \{p \colon A \to B \mid p \text{ is a fibration between fibrant objects in } \Leandspace_E \}\end{equation*}
 		and
 		\begin{equation*}\mathcal{Q} = \{q \colon A \to B \mid q \text{ is a trivial fibration between fibrant objects in } \Leandspace_E \},\end{equation*}
 		respectively.
 	\end{proposition}
 
 	As mentioned above, we may alternatively consider the slice model structure on $\dProfS/E$. It turns out that both model structures agree.
 	
 	\begin{lemma}\label{lemma:Both-sliced-model-structures-agree}
 		The model structure of \cref{proposition:Completed-model-structure-on-pro-slice} agrees with the slice model structure on $\dProfS/E$ under the equivalence $\Pro(\Leandspace_E) \simeq \dProfS/E$.
 	\end{lemma}
 
	\begin{proof}
		We will show that the (trivial) fibrations of both model structures agree. Note that $\dProfS/E$ is fibrantly generated by maps of the form $A \times E \to B \times E$, where $A \to B$ is a generating (trivial) fibration of $\dProfS$ (cf.\ \cref{equation:definition-Q-dSpace} and \cref{equation:definition-P-dSpace}). Note that if $A \to B$ is a generating (trivial) fibration of $\dProfS$, then it is a (trivial) fibration between fibrant objects in $\dSpace$. Since $E \times A \to E \times B$ lies in $\Leandspace_E$, it is a generating (trivial) fibration for the model structure of \cref{proposition:Completed-model-structure-on-pro-slice}.
		
		Conversely, suppose that $p \colon A \to B$ is a generating (trivial) fibration in the model structure of \cref{proposition:Completed-model-structure-on-pro-slice}. For $n \geq 0$, let $\pi_n \colon E \to \cosk_n E$ denote the projection and let $\pi^*_n \colon \Leandspace/\cosk_nE \to \Leandspace_E$ be the pullback functor. We leave it to the reader to verify that since $A$ and $B$ are objects of $\Leandspace_E$, there exists an $n \geq 0$ and a map $p' \colon A' \to B'$ in $\Leandspace/\cosk_nE$ such that $p$ is isomorphic to $\pi^*_np'$. We then have a diagram
		\begin{equation*}\begin{tikzcd}
			A \ar[r, two heads,"p"] \ar[dr, very near start,phantom,"\lrcorner"] \ar[d,two heads, "\sim" rot90] & B \ar[d,two heads, "\sim" rot90] \ar[dr, very near start,phantom,"\lrcorner"] \ar[r] & E \ar[d,two heads, "\sim" rot90, "\pi_n"] \\
			A' \ar[r,"p'"] & B' \ar[r] & \cosk_n E
		\end{tikzcd}\end{equation*}
		of pullback squares in $\dSpace$. The map $\pi_n$ is a trivial fibration in $\dSpace$ by the same argument as \cite[Lem.\ 5.3]{BlomMoerdijk2022ProfinitePartI}, hence all vertical maps in this diagram are trivial fibrations. Note that in the Reedy model structure on $\dSpace$, the generating trivial cofibrations can be chosen to have cofibrant codomains. By \cite[Prop.\ 4.3]{Hovey2004HomotopyComodules}, this holds for any of the model structures considered in \cref{ssec:Model-structures-preliminaries}. In particular, we may apply \cref{lemma:fibration-descent} below to deduce that $p' \colon A' \to B'$ is a (trivial) fibration between fibrant objects in $\gd\Space$. Since $A'$ and $B'$ are lean, this implies that $p'$ is a generating (trivial) fibration in $\dProfS$. In particular, $p'$ is a (trivial) fibration in $\dProfS/E$.
	\end{proof}

	\begin{lemma}\label{lemma:fibration-descent}
		Let $\mathcal{M}$ be a cofibrantly generated model category and suppose that there exists a set of generating trivial cofibrations with cofibrant codomain. Then for any pullback square
		\begin{equation*}\begin{tikzcd}
			M \ar[d,two heads,"p"'] \ar[r,"q'"] \ar[dr, very near start,phantom,"\lrcorner"] & M' \ar[d,"p'"] \\
			N \ar[r,two heads, "\sim","q"'] & N'
		\end{tikzcd}\end{equation*}
		where $p$ is a fibration and $q$ is a trivial fibration, the map $p'$ is a fibration.
	\end{lemma}

	\begin{proof}
		Let $X \trivcofarrow Y$ be a trivial cofibration with cofibrant codomain and suppose a lifting problem of the form
		\begin{equation*}\begin{tikzcd}
			X \ar[d, tail, "\sim" rot90] \ar[r] & M' \ar[d,"p'"] & M \ar[l,"q'"'] \ar[d,"p", two heads], \ar[dl,very near start,phantom,"\llcorner"] \\
			Y \ar[r,"f"'] \ar[ur,dashed] & N' & N \ar[l, two heads, "\sim"', "q"]
		\end{tikzcd}\end{equation*}
		is given. Since $Y$ is cofibrant, we may lift $f$ along $q$ to obtain a map $g \colon Y \to N$. The universal property of the pullback now provides a map $h \colon X \to M$. Letting $l$ be a lift in the diagram
		\begin{equation*}\begin{tikzcd}
			X \ar[d, tail, "\sim" rot90] \ar[r,"h"] & M \ar[d,"p"] \\
			Y \ar[r,"g"'] \ar[ur,dashed,"l"] & N
		\end{tikzcd}\end{equation*}
		the composite $q' l \colon Y \to M'$ is the desired lift.
	\end{proof}
 	
 	Let $\Leandspace_E^f$ denote the full simplicial subcategory $\Leandspace_E$ spanned by those objects that are fibrant in $\dSpace/E$, and let $N(\Leandspace_E^f)$ denote the homotopy-coherent nerve of this category. We write $\Pro_\infty(\Cin)$ for the $\infty$-categorical pro-category of an $\infty$-category $\Cin$, as defined in \cite[Def.\ 5.3.5.1]{Lurie2009HTT}. By \cite[Thm.\ A.7]{BlomMoerdijk2020SimplicialProPublished}, we may identify the underlying $\infty$-category of $\Pro(\Leandspace_E)$ as follows. 
 	
 	\begin{proposition}\label{proposition:Underlying-infty-category-slice}
 		The underlying $\infty$-category of $\dProfS/E$, and hence that of $\dProfS$, is equivalent to $\Pro_\infty(N(\Leandspace_E^f))$.
 	\end{proposition}
 
 	\begin{proof}
 		This follows by combining \cite[Thm.\ A.7]{BlomMoerdijk2022ProfinitePartI} with \cref{lemma:Both-sliced-model-structures-agree} and the Quillen equivalence between $\dProfS/E$ and $\dProfS$.
 	\end{proof}
 
	We will now identify the $\infty$-category $N(\Leandspace_E^f)$ in the case that $\dSpace$ is endowed with the model structure for complete dendroidal Segal spaces. Combined with the previous propositions, this yields the desired description of the underlying $\infty$-categories of $\dProfS_{RSC}$ and $\dSet$.
	
	\begin{definition}
		Let $\Op_\infty$ denote the $\infty$-category of $\infty$-operads.
		\begin{enumerate}[(a),noitemsep]
			\item The full $\infty$-subcategory $\Op_\infty^\pi$ of $\Op_\infty$ consists of those $\infty$-operads that are $\pi$-finite in the sense of \cite[Def.\ 2.16]{BlomMoerdijk2022ProfinitePartI} (cf.\ \cref{definition:pi-finite-complete-Segal-space}).
			\item The $\infty$-category $\Op_\infty^{\pi,\mathrm{o}}$ is defined as the full $\infty$-subcategory of $\Op_\infty^{\pi}$ spanned by those $\infty$-operads that moreover have no nullary operations.
			\item The full $\infty$-subcategory $\Op_\infty^{\pi,\mathrm{or}} \subset \Op_\infty^{\pi,\mathrm{o}}$ consists of those $\infty$-operads in which moreover all $1$-ary operations are invertible.
			\item The $\infty$-category $\Op_\infty^{\pi, \mathrm{cl}}$ is defined as the full subcategory of $\Op_\infty$ spanned by those $\infty$-operads that are both closed and $\pi$-finite in the sense of \cite[Def.\ 2.28]{BlomMoerdijk2022ProfinitePartI} (cf.\ \cref{definition:pi-finite-complete-Segal-space}).
		\end{enumerate}
	\end{definition}

	\begin{lemma}\label{lemma:Identifying-nerve-leandspaces}
		Let $\Leandspace_E^f$ denote the full simplicial subcategory of $\Leandspace$ spanned by those objects that are fibrant in $\dSpace_{RSC}$. Then the homotopy coherent nerve $N(\Leandspace_E^f)$ is equivalent to
		\begin{itemize}
			\item $\Op^\pi_\infty$ in case $\OmegaA = \Omega$,
			\item $\Op_\infty^{\pi,\mathrm{o}}$ in case $\OmegaA = \Omegao$,
			\item $\Op_\infty^{\pi,\mathrm{or}}$ in case $\OmegaA = \Omegaor$,
			\item $\Op_\infty^{\pi,\mathrm{cl}}$ in case $\OmegaA = \Omegacl$.
		\end{itemize}
	\end{lemma}

	\begin{proof}
		We will treat the case $\OmegaA = \OmegaA$, the other cases are similar. The $\infty$-category $\Op_\infty$ can be obtained as the homotopy coherent nerve of the full simplicial subcategory $\dd\Space_{RSC}^{\mathrm{fc}}$ of $\dd\Space_{RSC}$ spanned by the fibrant-cofibrant objects. Let $\dd\Space_{\pi}^{\mathrm{fc}}$ denote the full subcategory of $\dd\Space_{RSC}^\mathrm{fc}$ spanned by the complete dendroidal Segal spaces that are $\pi$-finite in the sense of \cref{definition:pi-finite-complete-Segal-space}. Then $\Op_\infty^\pi$ is the homotopy-coherent nerve of $\dd\Space^\mathrm{fc}_{\pi}$.
		
		Now note that the total space of any object of $\Leandspace_E$ is $\pi$-finite. Namely, if $X \to E$ is an object of $\Leandspace_E^f$, then by \cref{example:Writing-slice-over-E-as-pro-category} we may write $X$ as the pullback $X = p^*_nY$ of some map $Y \to \cosk_n E$ between lean dendroidal spaces. By \cref{lemma:fibration-descent}, the map $Y \to \cosk_nE$ is a fibration, hence $Y$ is fibrant and lean. By \cref{theorem:Characterization-lean-complete-dendroidal-Segal-spaces}, $Y$ is $\pi$-finite, so since $X = Y \times_{\cosk_n E} E \trivfibarrow Y$ is a trivial fibration, $X$ is $\pi$-finite as well. In particular, we have a forgetful functor
		\begin{equation*}U \colon \Leandspace_E^f \to \dd\Space^\mathrm{fc}_\pi.\end{equation*}
		The result follows if we show that this functor is a Dwyer--Kan equivalence of simplicial categories. To see that it is homotopically essentially surjective, let $X$ in $\dd\Space^{\mathrm{fc}}_\pi$ be given. By \cref{theorem:Characterization-lean-complete-dendroidal-Segal-spaces}, we may assume that $X$ is lean. Then $X \times E$ is an object of $\Leandspace_E^f$ and $X \times E \simeq X$, hence $U$ is homotopically essentially surjective. On the other hand, this map is homotopically fully faithful since $E \trivfibarrow *$ is a trivial fibration.
	\end{proof}
	
	\begin{theorem}\label{theorem:Underlying-infty-categories}
		\begin{enumerate}[(a),noitemsep]
			\item The underlying $\infty$-categories of $\dd\ProfS_{RSC}$ and $\dd\wh\Set$ are equivalent to $\Pro_\infty(\Op_\infty^\pi)$.
			\item The underlying $\infty$-categories of $\od\ProfS_{RSC}$ and $\od\wh\Set$ are equivalent to $\Pro_\infty(\Op_\infty^{\pi,\mathrm{o}})$.
			\item The underlying $\infty$-categories of $\cd\ProfS_{RSC}$ and $\cd\wh\Set$ are equivalent to $\Pro_\infty(\Op_\infty^{\pi,\mathrm{cl}})$.
			\item The underlying $\infty$-category of $\rod\ProfS_{RSC}$ is equivalent to $\Pro_\infty(\Op_\infty^{\pi,\mathrm{or}})$.
		\end{enumerate}
	\end{theorem}

	\begin{proof}
		Combine \cref{theorem:Profinite-Segal-spaces-vs-profinite-operads,proposition:Underlying-infty-category-slice,lemma:Identifying-nerve-leandspaces}.
	\end{proof}

%% file: Sections/BoavidaHorelRobertson.tex
\section{Levelwise profinite completion of \texorpdfstring{$\infty$}{infinity}-operads}\label{sec:Boavida-Horel-Robertson}
	
	The aim of this section is to compare the profinite completion of an $\infty$-operad as defined in \cref{definition:Profinite-completion-operad} with the profinite completion as defined by Boavida--Horel--Robertson in \cite{BoavidaHorelRobertson2019GenusZero}. This comparison will follow from a classification of the fibrant objects $X$ in $\rod\Space_{RS}$ satisfying $X(\eta) \simeq *$.
	
	\subsection{Comparison of different notions of profinite completion for \texorpdfstring{$\infty$}{infinity}-operads}\label{ssec:representability-problem-profinite-completion-operad}
	
	Let us briefly recall the definition of the profinite completion of an operad given by Boavida--Horel--Robertson in \cite[\S 5]{BoavidaHorelRobertson2019GenusZero}. First, they define an (uncoloured) profinite $\infty$-operad as a dendroidal profinite space $X$ such that the map $X(\eta) \to *$ is a weak equivalence and which moreover satisfies the \emph{Segal condition}; i.e.\ for any tree $R = S \circ_e T$ where $e$ is an internal edge of $R$, the Segal map
	\begin{equation*} X(R) \to X(S) \times_{X(e)}^h X(T)\end{equation*}
	is a weak equivalence in Quick's model structure for profinite spaces. Note that since $X(\eta) \simeq *$, the codomain of the Segal map is equivalent to the homotopy product $X(S) \times^h X(T)$. Let us write $\Op_\infty(\ProfS)$ for the full subcategory of $\dd\ProfS$ spanned by these uncoloured profinite $\infty$-operads. It is observed by Boavida--Horel--Robertson that the category $\dd\ProfS$ of dendroidal profinite spaces can be equipped with a model structure in which the weak equivalences are the maps that are levelwise an equivalence in $\s\wh\Set_Q$; we will use the Reedy model structure for this.
	They then suggest defining the profinite completion of an uncoloured $\infty$-operad $X$ as an object in $\Op_\infty(\ProfS)$ which corepresents the functor
	\begin{equation}\label{equation:Definition-profinite-completion-BHR}
		\Op_{\infty}(\ProfS) \to \s\Set; \quad Y \mapsto \mathbb{R}\Map(X,UY_f)
	\end{equation}
	up to weak equivalence. Here $Y_f$ denotes a fibrant replacement of $Y$ in the Reedy model structure on $\dd\ProfS$ and $UY_f$ denotes its underlying dendroidal space. An object $Z$ in $\Op_\infty(\ProfS)$ is said to corepresent this functor up to weak equivalence if there exists a zigzag of natural weak equivalences between this functor and the functor $Y \mapsto \mathbb{R}\Map(Z,Y)$. Since $\mathbb{R}\Map(-,-)$ is computed with respect to the levelwise equivalences in $\dProfS$ and since $Y_f$ need not be fibrant in $\dProfS_{RSC}$, it does not immediately follow that the profinite completion defined in \cref{definition:Profinite-completion-operad} agrees with the profinite completion in the sense of \cite{BoavidaHorelRobertson2019GenusZero}, given that the latter exists. It turns out that in many cases they do agree, as we will show in \cref{theorem:Main-theorem-Boavida-Horel-Robertson} below.
	
	In specific examples, such as when $X$ is the operad of little 2-discs, Boavida--Horel--Robertson show that a profinite completion in their sense exists and that it can be constructed levelwise. To this end, let us call an uncoloured $\infty$-operad $X$ \emph{profinitely adapted} if it has the property that for any $n_1,\ldots,n_k \geq 0$, the map
	\begin{equation*}
		(X(n_1) \times \cdots \times X(n_k))^\wedge \to \wh{X(n_1)} \times^h \cdots \times^h \wh{X(n_k)}
	\end{equation*} 
	is a weak equivalence. Here $\wh{(\cdot)}$ and $(\cdot)^\wedge$ denote the profinite completion functor. For conditions on the spaces of operations that ensure $X$ is profinitely adapted, the reader is referred to \cite[Prop.\ 3.9]{BoavidaHorelRobertson2019GenusZero} and \cite{Haine2024ProfiniteCompletionsProducts}. If one models a profinitely adapted $\infty$-operad $X$ as a dendroidal Segal space with $X(\eta) \simeq *$ and applies the profinite completion functor levelwise, then one clearly obtains a dendroidal profinite space $\wh X$ that lies in $\Op_{\infty}(\ProfS)$. Since the levelwise profinite completion functor is (derived) left adjoint to the forgetful functor $U \colon \dd\ProfS \to \dd\Space$, it follows that $\wh X$ corepresents the functor \cref{equation:Definition-profinite-completion-BHR} up to weak equivalence.
	
	The goal of this section is to show that the profinite completion of an uncoloured $\infty$-operad in the sense of \cite{BoavidaHorelRobertson2019GenusZero} always exists and that it is given by the profinite completion we constructed in \cref{definition:Profinite-completion-operad}. In our proof, we will slightly change our indexing category and use reduced open trees, instead of all trees.\footnote{Our reason for working with reduced open dendroidal spaces is of a technical nature: our proof makes extensive use of a certain filtration of $\Omegaor$ by subcategories $\Omegaorf{n}$ that is not as well behaved when working with the category $\Omega$ of all trees (cf.\ \cref{remark:Reason-for-using-outer-Reedy}).} A reduced open dendroidal Segal space $X$ is called \emph{profinitely adapted} if $X(\eta) \simeq *$ and for every $n_1, \ldots, n_k \geq 0$, the map
	\begin{equation*}
		(X(C_{n_1}) \times \cdots \times X(C_{n_k}))^\wedge \to \wh{X(C_{n_1})} \times^h \cdots \times^h \wh{X(C_{n_k})}
	\end{equation*} 
	is an equivalence.
	
	\begin{theorem}\label{theorem:Main-theorem-Boavida-Horel-Robertson} Let $\Op_\infty(\ProfS)$ denote the full subcategory of $\rod\ProfS$ spanned by the objects $X$ for which $X(\eta) \simeq *$ and which moreover satisfy the Segal condition.
		\begin{enumerate}[(a)]
			\item\label{item1:Main-theorem-Boavida-Horel-Robertson} A reduced open dendroidal profinite space $X$ lies in the full subcategory $\Op_\infty(\ProfS)$ if and only if $X(\eta) \simeq *$ and it is levelwise equivalent to a fibrant object in $\rod\ProfS_{RS}$.
			\item\label{item2:Main-theorem-Boavida-Horel-Robertson} Suppose $X$ is a profinitely adapted reduced open dendroidal Segal space. Then its profinite completion in the sense of \cref{definition:Profinite-completion-operad} agrees up to weak equivalence with its profinite completion as defined by Boavida--Horel--Robertson \cite[\S 5]{BoavidaHorelRobertson2019GenusZero}.
			\item\label{item3:Main-theorem-Boavida-Horel-Robertson} More generally, for any reduced open dendroidal Segal space $X$ with $X(\eta) \simeq *$ and any $Y$ in $\Op_\infty(\ProfS)$ that is Reedy fibrant, there is an equivalence
			\begin{equation*}\mathbb{R}\Map(\wh X_f, Y) \simeq \mathbb{R}\Map(X, UY)\end{equation*}
			that is natural in $Y$, where $\wh X_f$ denotes a fibrant replacement in $\rod\ProfS_{RS}$ of the levelwise profinite completion $\wh X$ of $X$. Here $\mathbb{R}\Map$ denotes the derived mapping space computed with respect to the levelwise equivalences.
		\end{enumerate}
	\end{theorem}
	
	Note that item \ref{item3:Main-theorem-Boavida-Horel-Robertson} states that for any reduced open dendroidal Segal space $X$ with $X(\eta) \simeq *$, its profinite completion in the sense of \cite[\S 5]{BoavidaHorelRobertson2019GenusZero} exists and is given by \cref{definition:Profinite-completion-operad}.
	
	\begin{remark}
		Since the square
		\begin{equation*}\begin{tikzcd}
			\rod\Space \ar[r,"\wh{(\cdot)}"] \ar[d,hook] & \rod\ProfS \ar[d,hook]\\
			\dd\Space \ar[r,"\wh{(\cdot)}"] & \dd\ProfS
		\end{tikzcd}\end{equation*}
		commutes up to natural equivalence, it follows from \cref{theorem:Main-theorem-Boavida-Horel-Robertson} that for any complete dendroidal Segal space $X$ with the property that $X(C_1) \simeq X(\eta) \simeq *$ and $X(\ol \eta) = \varnothing$, its profinite completion in the sense of \cref{definition:Profinite-completion-operad} must agree with that of \cite[\S 5]{BoavidaHorelRobertson2019GenusZero}. Combined with \cref{remark:Profinite-completions-agree}, it follows that these profinite completions can also be computed at the level of dendroidal (profinite) sets.
	\end{remark}
	
	We will deduce the theorem from the following proposition.
	
	\begin{proposition}\label{proposition:Main-proposition-Boavida-Horel-Robertson}
		Let $X \in \rod\ProfS$ be a Reedy fibrant reduced open dendroidal profinite space with $X(\eta) \simeq *$ and suppose that $X$ satisfies the Segal condition. Then $X$ is fibrant in $\rod\ProfS_{RS}$.
	\end{proposition}

	\begin{remark}
		Note that from the construction of the model structure $\rod\ProfS_{RS}$, it is not at all obvious that a Reedy fibrant object in $\rod\ProfS$ that satisfies the Segal condition is fibrant in $\rod\ProfS_{RS}$.
	\end{remark}
	
	\Cref{ssec:Boavida-Horel-Robertson-Main-Proof} below is devoted to the proof of this proposition, but let us first show how to deduce \cref{theorem:Main-theorem-Boavida-Horel-Robertson} from it.
	
	\begin{proof}[Proof of \cref{theorem:Main-theorem-Boavida-Horel-Robertson}, assuming \cref{proposition:Main-proposition-Boavida-Horel-Robertson}]
		 For the ``if'' direction of \ref{item1:Main-theorem-Boavida-Horel-Robertson}, note that by \cref{theorem:Profinite-spaces-equivalence-detected-underlying}, a Reedy fibrant object $X$ in $\rod\ProfS$ satisfies the Segal condition if and only if the underlying reduced open dendroidal space $UX$ does. Since $U \colon \rod\ProfS_{RS} \to \rod\Space_{RS}$ is right Quillen, any fibrant object in $\rod\ProfS$ satisfies the Segal condition. For the ``only if'' direction, suppose $X$ lies in $\Op_\infty(\ProfS)$. Then $X$ is levelwise equivalent to a Reedy fibrant object, which is fibrant in $\rod\ProfS_{RS}$ by \cref{proposition:Main-proposition-Boavida-Horel-Robertson}.
		
		To deduce item \ref{item2:Main-theorem-Boavida-Horel-Robertson}, let $X$ be a profinitely adapted reduced open dendroidal Segal space. Then its levelwise profinite completion $\wt X$ satisfies the Segal condition, hence it lies in $\Op_\infty(\ProfS)$. By \ref{item1:Main-theorem-Boavida-Horel-Robertson}, $\wt X$ is levelwise equivalent to an object $Y$ that is fibrant in $\rod\ProfS_{RS}$. Since this object is a fibrant replacement of $\wt X$ in $\rod\ProfS_{RS}$, we conclude that the profinite completion of $X$ as defined in \cref{definition:Profinite-completion-operad} is levelwise equivalent to that of \cite[\S 5]{BoavidaHorelRobertson2019GenusZero}.
		
		For item \ref{item3:Main-theorem-Boavida-Horel-Robertson}, it suffices to show that for any normal $X$ and any Reedy fibrant $Y$ in $\Op_\infty(\ProfS)$, the natural map
		\begin{equation*}\sHomKan(\wh X_f, Y) \to \sHomKan(\wh X, Y) \cong \sHomKan(X, UY) \end{equation*}
		is a weak equivalence, where $\wh X \wearrow \wh X_f$ is a fibrant replacement of $\wh X$ in $\rod\ProfS_{RS}$. This follows since $\wh X \to \wh X_f$ is a weak equivalence in $\rod\ProfS_{RS}$ and $Y$ is fibrant in $\rod\ProfS_{RS}$ by \cref{proposition:Main-proposition-Boavida-Horel-Robertson}.
	\end{proof}

	\subsection{A classification of fibrant objects in \texorpdfstring{$\rod\ProfS_{RS}$}{rodPSpaces\_\{RS\}}}\label{ssec:Boavida-Horel-Robertson-Main-Proof}
	
	Recall the filtration
	\begin{equation*}
		\{\eta\} = \Omegaorf{1} \subset  \Omegaorf{2} \subset \cdots \subset \cup_n \Omegaorf{n} = \Omegaor
	\end{equation*}
	from \cref{eq:filtration-Omegaor-by-arity} and write $w_n$ for the inclusion $\Omegaorf{n} \hookrightarrow \Omegaor$. This filtration yields for any reduced open dendroidal profinite space $X$ an isomorphism $X \cong \lim_n w_{n*}w_n^*X$, where $w^*_n$ and $w_{n*}$ denote the restriction and right Kan extension functors, respectively. Our strategy for proving \cref{proposition:Main-proposition-Boavida-Horel-Robertson} is to show inductively that $w_{n*}w_n^*X$ is fibrant and then conclude from this that $X$ must itself be fibrant. Unfortunately, this does not exactly work, since $w_{n*}w_n^*X$ need not be Reedy fibrant. This issue is solved by working with the outer Reedy model structure $\rod\ProfS_{oR}$ defined in \cref{ssec:Outer-Reedy-structure}. We will therefore prove the analogue of \cref{proposition:Main-proposition-Boavida-Horel-Robertson} for the outer Reedy structure, from which \cref{proposition:Main-proposition-Boavida-Horel-Robertson} is easily deduced. Recall that an object in $\rod\ProfS$ is called outer fibrant if it is fibrant in $\rod\ProfS_{oR}$.
	
	\begin{proposition}\label{proposition:Main-proposition-Boavida-Horel-Robertson-Outer-version}
		Let $X \in \rod\ProfS$ be an outer fibrant reduced open dendroidal profinite space with $X(\eta) \simeq *$ and suppose that $X$ satisfies the Segal condition. Then $X$ is fibrant in $\rod\ProfS_{oRS}$.
	\end{proposition}

	Recall from \cref{proposition:Existence-outer-Reedy-Segal-model-structure} that the weak equivalences of $\rod\ProfS_{oRS}$ agree with those of $\rod\ProfS_{RS}$. Since the outer Reedy structure on $\Omegaor$ has fewer positive morphisms than the standard Reedy structure, the class of fibrant objects in $\rod\ProfS_{oRS}$ contains the class of fibrant objects in $\rod\ProfS_{RS}$.

	\begin{proof}[Proof of \cref{proposition:Main-proposition-Boavida-Horel-Robertson}, assuming \cref{proposition:Main-proposition-Boavida-Horel-Robertson-Outer-version}]
		Let $X$ be a fibrant object in the standard Reedy model structure $\rod\ProfS_R$ which satisfies the Segal condition and suppose that $X(\eta) \simeq *$. Let $X \to \wt X$ be a fibrant replacement in $\rod\ProfS_{RS}$. Then $X \to \wt X$ is a weak equivalence between fibrant objects in $\rod\ProfS_{oRS}$ by \cref{proposition:Main-proposition-Boavida-Horel-Robertson-Outer-version}, hence a levelwise equivalence. By \cref{lemma:Detecting-fibrant-objects-in-left-Bousfield-localization}, it follows that $X$ is fibrant in $\rod\ProfS_{RS}$.
	\end{proof}

	The inductive proof that $w_{n*}w_n^*X$ is fibrant makes use of the intermediate full subcategory $\Omegaorfp{n}$ that was defined as the full subcategory $\Omegaor$ spanned by the objects of $\Omegaorf{n}$ and the $(n+1)$-corolla $C_{n+1}$. We have a diagram of inclusions
	\begin{equation*}
		\Omegaorf{n-1} \xhookrightarrow{u_n} \Omegaorfp{n-1} \xhookrightarrow{v_n} \Omegaorf{n} \xhookrightarrow{w_n} \Omegaor
	\end{equation*}
	inducing functors
	\begin{equation*}
		\begin{tikzcd}[column sep=scriptsize]
			\rod\ProfS \ar[r, shift left, "w_n^*"] & \rodProfSf{n} \ar[r, shift left, "v_n^*"] \ar[l, shift left, "w_{n*}"] & \rodProfSfp{n-1} \ar[r, shift left, "u_n^*"] \ar[l, shift left, "v_{n*}"] & \rodProfSf{n-1} \ar[l, shift left, "u_{n*}"]
		\end{tikzcd}
	\end{equation*}
	by restriction and right Kan extension. Recall from \cref{ssec:Outer-Reedy-structure} that $\rodProfSf{n}$ and $\rodProfSfp{n}$ denote the categories of $\s\wh\Set$-valued presheaves on $\Omegaorf{n}$ and $\Omegaorfp{n+1}$, respectively. We will show that these functors are right Quillen with respect to the outer Reedy model structures. To this end, let us make the following definition.
	
	\begin{definition}
		Let $\RReedy$ be a generalized Reedy category. A full subcategory $\Cin \subset \RReedy$ is called \emph{downwards closed} if for any $c \in \Cin$ and any $r \in \RReedy$,
		\begin{enumerate}[(1)]
			\item if there exists a positive morphism $r \to c$, then $r \in \Cin$, and
			\item if there exists a negative morphism $c \to r$, then $r \in \Cin$.
		\end{enumerate}
	\end{definition}
	
	\begin{remark}
		Suppose that $\RReedy$ is a Reedy category with degree function $d$. Then for any $m \in \bbN$, the full subcategory $\RReedy_{\leq m} \subset \RReedy$ spanned by those objects $r$ for which $d(r) \leq m$ is a downwards closed full subcategory. However, not every downwards closed full subcategory needs to be of this form. For example, $\Omegaorf{n} \subset \Omegaor$ is downwards closed with respect to the outer Reedy structure, but it cannot be of this form since there are only finitely many trees with weight at most $m$.
	\end{remark}
		
	\begin{lemma}\label{lemma:Restricting-along-downwards-closed-subcategory-latching-matching}
		Let $\RReedy$ be a generalized Reedy category and let $i \colon \Cin \hookrightarrow \RReedy$ be the inclusion of a downwards closed full subcategory. Then $i^*$ commutes with matching and latching objects. More precisely, if $\catE$ is a (co)complete category and $X \colon \RReedy \to \catE$ is any functor, then for any $c \in \Cin$, the maps
		\begin{equation*}
			L_c(i^*X) \to L_cX \quad \text{and} \quad M_cX \to M_c(i^*X)
		\end{equation*}
		are isomorphisms.
	\end{lemma}
	
	\begin{proof}
		The latching object $L_cX$ is computed as the colimit $\colim_{c \twoheadrightarrow r} X(r)$, indexed by the category of negative morphisms with source $c$ that are not isomorphisms. Since $\Cin$ is downwards closed, it follows that $r$ lies in $\Cin$. The statement about the matching object is dual.
	\end{proof}
	
	\begin{corollary}
		The functors $u_{n*}$, $v_{n*}$, $w_{n*}$, $u_n^*$, $v_n^*$ and $w^*_n$ are right Quillen with respect to the outer Reedy model structures.
	\end{corollary}
	
	\begin{proof}
		Since $\Omegaorf{n-1} \subset \Omegaorfp{n-1} \subset \Omegaorf{n} \subset \Omegaor$ is a sequence of inclusions of downwards closed full subcategories, the functors $u^*_n$, $v^*_n$ and $w^*_n$ commute with latching and matching objects by \cref{lemma:Restricting-along-downwards-closed-subcategory-latching-matching}. This implies that they preserve (trivial) outer Reedy (co)fibrations, hence that they are both left and right Quillen.
	\end{proof}

	\begin{remark}\label{remark:Reason-for-using-outer-Reedy}
	The functor $w_n^*$ is generally not right Quillen with respect to the standard Reedy structures. If this were the case, then $w_{n!}$ would preserve normal monomorphisms, which it does not. The outer Reedy model structure solves this issue by having fewer cofibrations.
	\end{remark}
	
	\begin{proposition}\label{lemma:Four-filtration-functors-are-right-Quillen}
		The functors $u_{n*}$, $u_n^*$, $v_n^*$ and $w^*_n$ are right Quillen with respect to the model categories of \cref{proposition:Existence-outer-Reedy-Segal-model-structure,proposition:Existence-outer-Reedy-Segal-model-structure-restricted}.
	\end{proposition}

	\begin{proof}
		By Corollary A.2 and Remark A.3 of \cite{Dugger2001Replacing}, it suffices to show that these functors preserve fibrant objects. In light of \cref{lemma:Detecting-fibrant-objects}, it suffices to prove that these functors preserve lean objects satisfying the Segal conditions. It is clear that $u_{n*}$ preserves lean objects. Moreover, note that $X \in \rodProfSfp{n-1}$ satisfies the Segal condition if and only if its restriction to $\rodProfSf{n-1}$ does, since the Segal condition for $C_n$ is vacuous. This shows that $u_{n*}$ also preserves objects satisfying the Segal condition.
		
		On the other hand, it is clear that $v_n^*$, $w_n^*$ and $u^*_n$ preserve objects satisfying the Segal condition, while it follows from \cref{lemma:Detecting-lean-using-outer-matching-maps} that they preserve lean objects.
	\end{proof}

	The functors $v_{n*}$ and $w_{n*}$ are also right Quillen with respect to these model structures, but this will require a careful analysis of these functors. Given a tree $T$ in $\Omegaor$, the value of $w_{n*}X$ at $T$ is given by the limit
	\begin{equation*}(w_{n*}X)(T) = \lim_{S \to T} X(S)\end{equation*}
	over the category $(\Omegaorf{n}/T)^{op}$ of all maps $S \to T$ with $S \in \Omegaorf{n}$. Given such a map, there is always a unique maximal subtree $R$ of $T$ with the property that $R$ lies in $\Omegaorf{n}$ and that $S \to T$ factors through $R \hookrightarrow T$. Namely, cut the tree $T$ at each inner edge $e$ that is attached to a vertex with more than $n$ ingoing edges to obtain a collection $\mathcal{T}$ of subtrees of $T$. Throwing away all trees in $\mathcal{T}$ that have a vertex with more than $n$ ingoing edges, what remains are exactly the maximal subtrees of $T$ that lie in $\Omegaorf{n}$. Since none of these subtrees share any edge, it follows that they span a cofinal discrete subcategory of $(\Omegaorf{n}/T)^{op}$. Denoting the set of these subtrees by $\mathrm{Max}_{(n)}(T)$, we see that
	\begin{equation*}(w_{n*}X)(T) = \prod_{R \in \mathrm{Max}_{(n)}(T)} X(R). \end{equation*}
	For $v_{n*}$ one can argue similarly, though the cofinal subcategory one obtains is not discrete. Given a tree $T$ in $\Omegaorf{n}$ with root edge $r$, let $E$ denote the set of all inner edges attached to a vertex with $n$ inputs. By cutting $T$ into smaller subtrees at these edges, we obtain a decomposition of $T$ into subtrees. For any $e \in E \cup \{r\}$, we let $R_e$ denote the subtree in this decomposition whose root edge is given by $e$. Now note that any map $S \to T$ where $S$ lies in $\Omegaorfp{n-1}$ factors through a unique subtree of the form $R_e$, unless $S = e$ for some $e \in E$, in which case it factors through precisely two such subtrees. By a similar argument as above, we obtain that $(v_{n*}X)(T)$ is an iterated pullback of the (profinite) spaces $X(R_e)$ over $X(e)$, which we denote as
	\begin{equation*}(v_{n*}X)(T) = \tilde{\bigtimes\limits_{e \in E \cup \{r\}}} X(R_e).\end{equation*}
	Using that pullbacks of the form $X(S) \times_{X(e)} X(T)$ are automatically homotopy pullbacks if $X$ is outer fibrant, we will deduce from these formulas that $v_{n*}$ and $w_{n*}$ preserve outer fibrant objects satisfying the Segal condition.
	
	\begin{proposition}\label{lemma:Other-2-filtration-functors-are-right-Quillen}
		The functors
		\begin{align*}
			&v_{n*}\colon \rodProfSfp{n-1}_{oRS} \to \rodProfSf{n}_{oRS}, \\
			&w_{n*} \colon \rodProfSf{n}_{oRS} \to \rod\ProfS_{oRS}
		\end{align*}
		are right Quillen.
	\end{proposition}

	\begin{proof}
		By the same argument as in \cref{lemma:Four-filtration-functors-are-right-Quillen}, it suffices to prove that $v_{n*}$ and $w_{n*}$ preserve outer fibrant objects satisfying the Segal condition. To this end, let $X$ be an object of $\rodProfSf{n}$ that is outer fibrant and satisfies the Segal condition, and suppose we are given a tree $R = S \circ_e T$ with $e$ an inner edge. Let $U \in \mathrm{Max}_{(n)}(R)$ be such that $e \in U$ and write $U = U_1 \circ_e U_2$. Note that $e$ need not be an inner edge of $U$, so it could be that $U_1 = e$ or $U_2 = e$. The map
		\begin{equation*}
			(w_{n*}X)(R) \to (w_{n*}X)(S) \times_{(w_{n*}X)(e)} (w_{n*}X)(T)
		\end{equation*}
		can be identified with the map
		\begin{equation*}
			\prod_{M \in \mathrm{Max}_{(n)}(R)}  X(M) \to \left(X(U_1) \times_{X(e)} X(U_2)\right) \times \prod_{\substack{M \in \mathrm{Max}_{(n)}(R) \\ M \neq U}} X(M)
		\end{equation*}
		which is an isomorphism on every factor except possibly for the one corresponding to $U$, where it is the Segal map $X(U) \to X(U_1) \times_{X(e)} X(U_2)$. If $e$ is not an inner edge in $U$ then this map is an isomorphism, and otherwise it is a weak equivalence since $X$ satisfies the Segal condition. We conclude that $w_{n*}X$ satisfies the Segal condition.
		
		The argument for $v_{n*}$ is similar, replacing the product with an iterated pullback.
	\end{proof}

	\begin{remark}\label{remark:vn*-is-a-right-Quillen-equivalence}
		In fact, the functor $v_{n*}$ is a right Quillen equivalence. Since $v_n^*$ is also right Quillen, it suffices to show that the counit $v^*_nv_{n*} X \to X$ is a weak equivalence for any fibrant $X$ in $\rodProfSfp{n-1}$ and that the unit $Y \to v_{n*} v_n^*Y$ is a weak equivalence for any fibrant $Y$ in $\rodProfSf{n}$. The first follows since the counit is an isomorphism. For the unit, note that since $Y$ satisfies the Segal condition, so does $v_n^*Y$. By \cref{lemma:Other-2-filtration-functors-are-right-Quillen}, we see that $v_{n*}v_n^*Y$ satisfies the Segal condition. In particular, $Y \to v_{n*}v_n^*Y$ is a map between reduced open dendroidal spaces satisfying the Segal condition, hence it is a weak equivalence if and only if $Y(C_m) \to (v_{n*}v_n^*Y)(C_m)$ is a weak equivalence for every $m \leq n$. This is the case since these maps are isomorphisms.
	\end{remark}

	Recall that our aim is to prove that $w_{n*}w_n^*X$ is fibrant for every $n$ and then use \cref{lemma:Detecting-fibrant-objects} to conclude that $X = \lim_n w_{n*}w_n^*$ is fibrant. Since $w_{n*}$ is right Quillen, it suffices to show that $w_n^*X$ is fibrant in $\rodProfSf{n}_{oRS}$ for every $n$. This will be done inductively, using \cref{corollary:Reedy-extension-fact} to extend an object of $\rodProfSf{n-1}$ to one of $\rodProfSfp{n-1}$ in a very controlled way.
	
	Given a functor $F \colon \Cin \to \Din$ between categories and an object $d \in \Din$, we write $F^{-1}(d)$ for the \emph{fiber category} of $F$ above $d$; that is, $F^{-1}(d)$ is the (non-full) subcategory of $\Cin$ whose objects are the $c \in \Cin$ such that $F(c) = d$, and whose maps are those $f \colon c \to c'$ for which $F(f) = \id_d$.
	
	\begin{lemma}\label{corollary:Reedy-extension-fact}
		Let $X \in \rodProfSf{n-1}$ with $X(\eta) \cong *$ be given. Then the fiber category $\mathbf{E} \coloneqq (u^*_n)^{-1}(X)$ is equivalent to the slice category $\Lo_{C_n} X/(\eqswhSet{\Sigma_n})$, where the equivalence $\mathbf{E} \to \Lo_{C_n} X/(\eqswhSet{\Sigma_n})$ is given by sending $Y \in \rodProfSfp{n-1}$ to the map $\Lo_{C_n} X = \Lo_{C_n} Y \to Y(C_n)$.
	\end{lemma}
	
	\begin{proof}
		By essentially the same argument as \cite[Lem.\ VII.2.9]{GoerssJardine2009SimplicialHomotopyTheory}, we see that the category $\mathbf{E}$ is equivalent to the category of factorizations of the latching-to-matching map $\Lo_{C_n} X \to \Mo_{C_n} X$. Since $\Mo_{C_n} X = X(\eta)^{n+1} \cong *$, the result follows.
	\end{proof}

	The following result is the core ingredient of the proof of \cref{proposition:Main-proposition-Boavida-Horel-Robertson-Outer-version}.

	\begin{proposition}\label{proposition:The-induction-step}
		Let $Y \in \rodProfSfp{n-1}$ with $Y(\eta) \cong *$ be given such that $u^*_nY$ is fibrant in $\rodProfSf{n-1}_{oRS}$ and $Y(C_n)$ is fibrant in Quick's model structure on $\eqswhSet{\Sigma_n}$. Then $Y$ is fibrant in $\rodProfSfp{n-1}_{oRS}$.
	\end{proposition}
	
	\begin{proof}
		By \cref{proposition:Existence-outer-Reedy-Segal-model-structure-restricted}, we may write $u^*_nY$ as a cofiltered limit $\lim_{i \in I} X_i$ of lean outer fibrant objects satisfying the Segal condition. Moreover, by \cref{lemma:Profinite-G-Space-Filtered-colimit-lean-G-Kan} we may write $Y(C_n)$ as a cofiltered limit $\lim_{j \in J} Z_j$ of lean $\Sigma_n$-Kan complexes. Since the latching objects for the outer Reedy structure on $\Omegaor$ are computed as finite colimits, we see that $\Lo_{C_n} Y = \Lo_{C_n} (\lim_i X_i) \cong \lim_i \Lo_{C_n} X_i$. Now consider the map
		\begin{equation*}
			\lim_i \Lo_{C_n} X_i \cong \Lo_{C_n} Y \to Y(C_n) = \lim_j Z_j.
		\end{equation*}
		The construction of Proposition 3.1 of \cite[Appendix]{ArtinMazur1986Etale} yields a cofiltered category $K$ together with cofinal functors $\theta \colon K \to I$ and $\nu \colon K \to J$ and natural maps $\Lo_{C_n} X_{\theta(k)} \to Z_{\nu(k)}$ such that
		\begin{equation*}
			\lim_{k \in K} (\Lo_{C_n} X_{\theta(k)} \to Z_{\nu(k)}) = (\Lo_{C_n} Y \to Y(C_n)).
		\end{equation*}
		By applying \cref{corollary:Reedy-extension-fact} to the maps $\Lo_{C_n} X_{\theta(k)} \to Z_{\nu(k)}$, we obtain a diagram $\{W_k\}_{k \in K}$ in $\rodProfSfp{n-1}$ such that $u^* W_k = X_{\theta(k)}$ and $W_k(C_n) = Z_{\nu(k)}$. In particular, we see that $\lim_k W_k \cong Y$. Moreover, $W_k$ is lean by \cref{lemma:Detecting-lean-using-outer-matching-maps} since $X_{\theta(k)}$ and $Z_{\nu(k)}$ are. The object $W_k$ is outer fibrant since $X_{\theta(k)}$ is outer fibrant and $Z_{\nu(k)}$ is fibrant in $\eqswhSet{\Sigma_n}$, while it satisfies the Segal condition since $X_{\theta(k)}$ does. We conclude by \cref{lemma:Detecting-fibrant-objects} that $Y \cong \lim_k W_k$ is fibrant.
	\end{proof}

	Let us now paste everything together and deduce \cref{proposition:Main-proposition-Boavida-Horel-Robertson-Outer-version}, completing the proof of \cref{theorem:Main-theorem-Boavida-Horel-Robertson}.
	
	\begin{proof}[Proof of \cref{proposition:Main-proposition-Boavida-Horel-Robertson-Outer-version}]
		Let an outer fibrant object $X$ in $\rod\ProfS$ be given and suppose that $X(\eta) \simeq *$ and that $X$ satisfies the Segal condition. Choose an element $x \in X(\eta)$ and define $X'$ via the pullback
		\begin{equation*}\begin{tikzcd}
			X'(T) \ar[r," "] \ar[d," "] \ar[dr,"\lrcorner",very near start,phantom] & X(T) \ar[d," ",two heads]\\
			* \ar[r,"{(x_0,\ldots,x_0)}","\sim"'] & \prod_{e \in E(T)} X(\eta),
		\end{tikzcd}\end{equation*}
		where $E(T)$ denotes the set of edges of $T$.
		Then $X' \to X$ is a levelwise equivalence, $X'(\eta) \cong *$ and $X'$ is again outer fibrant. By \cref{lemma:Detecting-fibrant-objects-in-left-Bousfield-localization}, it suffices to show that $X'$ is fibrant in $\rod\ProfS_{oRS}$, so we assume without loss of generality that $X(\eta) \cong *$.
		
		Since $X = \lim_n w_{n*}w_n^*X$, it suffices by \cref{lemma:Detecting-fibrant-objects} and \cref{lemma:Other-2-filtration-functors-are-right-Quillen} to show that $w_n^*X$ is fibrant in $\rodProfSf{n}$ for every $n \geq 1$. It is clear that $w_1^*X$ is fibrant in $\rodProfSf{1}_{oRS}$ since $X(\eta) \cong *$ and $\Omegaorf{1} = \{\eta\}$. Now suppose that $w^*_{n-1}X$ is fibrant in $\rodProfSf{n-1}$. Since $X$ is outer fibrant, we see that $X(C_n)$ is fibrant in Quick's model structure on $\eqswhSet{\Sigma_n}$. By \cref{proposition:The-induction-step}, $v^*_nw^*_nX$ is fibrant in $\rodProfSfp{n-1}$ and by the argument of \cref{remark:vn*-is-a-right-Quillen-equivalence}, the unit $w_n^*X \to v_{n*}v_n^*w_n^*X$ is a levelwise weak equivalence. Since $v_{n*}v_n^*w_n^*X$ is fibrant by \cref{lemma:Other-2-filtration-functors-are-right-Quillen}, it follows from \cref{lemma:Detecting-fibrant-objects-in-left-Bousfield-localization} that $w_n^*X$ is fibrant in $\rodProfSf{n}_{oRS}$.
	\end{proof}